\documentclass[11pt,A4paper]{article}
\usepackage{latexsym,amsfonts,amsmath,amssymb,graphicx,theorem}
\newtheorem{theo}{Theorem}[section]
\newtheorem{defi}[theo]{Definition}
\newtheorem{lem}[theo]{Lemma}

\newtheorem{prop}[theo]{Proposition}

\newtheorem{rem}[theo]{Remark}
\begin{document}

\begin{center}
\textbf{{\LARGE Spectral types of linear $q$-difference equations and $q$-analog of middle convolution}}
\end{center}

\begin{center}
{\large Hidetaka Sakai and Masashi Yamaguchi

Graduate School of Mathematical Sciences, The university of Tokyo, \\
Komaba, Tokyo 153-8914, Japan.}
\end{center}
\begin{abstract}
We give a $q$-analog of middle convolution for linear $q$-difference equations with
rational coefficients. In the differential case, middle convolution is defined by Katz, and he examined properties of middle convolution in detail. In this paper, we define a $q$-analog of middle convolution. Moreover, we show that it also can be expressed as a $q$-analog of Euler transformation. The $q$-middle convolution transforms Fuchsian type equation to Fuchsian type equation and preserves rigidity index of $q$-difference equations.\\
\end{abstract}


\section{Introduction}
{\ }

In this paper, we give a $q$-analog of middle convolution for linear $q$-difference equations with
rational coefficients, and we show properties of the $q$-middle convolution. Before that, we briefly look over the theory of middle convolution for differential equations.

At first, we look over a theory of Katz in \cite{K}. He defined addition and middle convolution for solutions of differential equations of Schlesinger normal form
\begin{equation}
\frac{dY}{dx}(x)=A(x)Y(x),\ \ A(x)=\sum _{k=1}^N\frac{A_k}{x-t_k}\ \ (t_k\in \mathbb{C},\ A_k\in \mathrm{M}_m(\mathbb{C})).
\end{equation}

These operations transform Fuchsian equation to Fuchsian equation and preserve rigidity index of the equation. Rigidity index is the integer related to the number of accessory parameters. Accessory parameters are parameters which are independent of eigenvalues of $A_k,\, A_{\infty}=-(A_1+\cdots +A_N))$. If the equation (1) has no accessory parameters, it is called ``rigid". Katz showed that any irreducible rigid Fuchsian equations can be obtained from a certain 1st order equation by finite iterations of additions and middle convolutions. Katz's theorem tells that there exists integral representation of solutions of any irreducible rigid Fuchsian equations, because an addition transforms solution $Y(x)$ of the equation (1) to
\[ \prod _{k=1}^r(x-a_k)^{b_k}\cdot Y(x)\, (a_k,b_k\in \mathbb{C})\]
and a middle convolution is integral transformation for solution $Y(x)$ of the equation (1). 

\begin{rem}\normalfont{There are two types}, ``additive version" and ``multiplicative version" of middle convolution defined by Katz. Additive version is transformation for equations. Multiplicative version is transformation for solutions. Multiplicative middle convolution induces a transformation of monodromy representation. In this paper, we treat the similar version to the former, which should be called ``additive version" $q$-middle convolution. In the $q$-difference case, we think that connection matrix between two local solutions at singularities $x=0,\infty$ correspond to monodromy of differential equation. Birkhoff studied the connection matrix $P(x)$ for local solutions $Y_0(x),Y_{\infty}(x)$ at singularities $x=0,\infty$ of linear $q$-difference system with polynomial coefficient $Y(qx)=A(x)Y(x)$. Furthermore, Sauloy considered a category of linear $q$-difference systems with rational coefficients, a category of solutions and a category of connection data in \cite{Sa}. He gave Riemann-Hilbert correspondence for these categories. Based on the Sauloy's result, Roques studied rigidity of connection of linear $q$-difference systems with rational coefficients in \cite{Ro}.\ $\square$
\end{rem}

We referred to an easier construction of Dettweiler and Reiter in order to define the $q$-analog of middle convolution. Let us look over a result of Dettweiler and Reiter in \cite{DR1,DR2}. They express Katz's middle convolution in terms of matrices. The next transformation is called ``convolution" with parameter $\lambda \in \mathbb{C}$:
\begin{align}
&\frac{dZ}{dx}(x)=G(x)Z(x),\ \ G(x)=\sum _{k=1}^N\frac{G_k}{x-t_k}\ \ (G_k\in \mathrm{M}_{mN}(\mathbb{C})),\\
&G_k=\begin{pmatrix} & & O & & \\ A_1 & \dotsi & A_k+\lambda 1_m & \dotsi & A_N \\ & & O & & \end{pmatrix} (k\, \mathrm{th\ entry})\ \ (1 \le k \le N,\ 1_m=\{ \delta _{i,j}\}_{1\le i,j\le m}\in \mathrm{M}_m(\mathbb{C})).
\end{align}

Moreover, we define two linear spaces

\begin{equation}
\mathcal{K}=\begin{pmatrix} \ker A_1 \\
 \vdots \\
 \ker A_N  \end{pmatrix},\ \ \ \ \mathcal{L}=\ker (G_1+\cdots +G_N).
\end{equation}

Let $\overline{G}_k$ be a matrix induced by the action of $G_k$ on $\mathbb{C}^{mN}/(\mathcal{K}+\mathcal{L})$. We define middle convolution
\[ mc_{\lambda}\, :\ (A_1,\, \ldots ,\, A_n)\longmapsto (\overline{G}_1,\, \ldots ,\, \overline{G}_n).\]

We obtain a similar transformation by considering the Dettweiler and Reiter's setting in the $q$-difference case.

Let 
\begin{align*}
\boldsymbol{B}&={}^t(B_1, \ldots ,B_N,B_{\infty}) \in (\mathrm{M}_m(\mathbb{C}))^{N+1},\\
\boldsymbol{b}&={}^t(b_1, \ldots ,b_N) \in (\mathbb{C}\backslash \{ 0 \} )^N\ \ (b_i=b_j \Rightarrow i=j).
\end{align*}
We set an equation
\begin{equation}
E_{\boldsymbol{B},\boldsymbol{b}}\, :\ \sigma _x Y(x)=B(x)Y(x),\ \ \ B(x)=B_{\infty}+\sum _{i=1}^N \frac{B_i}{1-\frac{x}{b_i}}.
\end{equation}

For an equation $E_{\boldsymbol{B},\boldsymbol{b}}$, we define the $q$-convolution.

\begin{defi}\label{conv}$(q$-$\mathrm{convolution})$\label{conv}\ \ Let $\mathcal{E}$ be the set of $E_{\boldsymbol{B},\boldsymbol{b}}$'s. For $E_{\boldsymbol{B},\boldsymbol{b}} \in \mathcal{E},\lambda \in \mathbb{C}$, we define $q$-convolution$\, c_{\lambda}:\mathcal{E}\longrightarrow \mathcal{E}\ (E_{\boldsymbol{B},\boldsymbol{b}} \longmapsto E_{\boldsymbol{F},\boldsymbol{b}})$ as
\begin{equation}
\begin{split}
\boldsymbol{F}&=(F_1, \ldots ,F_N,F_{\infty}) \in (\mathrm{M}_{(N+1)m}(\mathbb{C}))^{N+1},\\
F_i&=\begin{pmatrix} & & O & & \\ B_0 & \dotsi & B_i-(1-q^{\lambda})1_m & \dotsi & B_N \\ & & O & & \end{pmatrix} (i+1\, \mathrm{th\ entry}),\ 1 \le i \le N,\\
F_{\infty}\! &=1_{(N+1)m}-\widehat{F},\\
\widehat{F}&=(B_{t-1})_{1 \leq s,t \leq N+1}=\begin{pmatrix} B_0 & \dotsi & B_N \\ \vdots & \ddots & \vdots \\ B_0 & \dotsi & B_N \end{pmatrix},\ B_0=1_m-B_{\infty}-\sum _{j=1}^NB_j.
\end{split}
\end{equation}
\end{defi}

Furthermore, we define the $q$-middle convolution.

\begin{defi}$(q$-$\mathrm{middle\ convolution})$\ \ Let $\mathcal{V}= \mathbb{C}^m$ and $\boldsymbol{F}$-invariant subspaces of $\mathcal{V}^{N+1}$ as
\begin{equation}
\mathcal{K}= \mathcal{K}_{\mathcal{V}}=\bigoplus _{i=0}^N\ker B_i,\ \mathcal{L}=\mathcal{L}_{\mathcal{V}}(\lambda )=\ker (\widehat{F}-(1-q^{\lambda})1_{(N+1)m}).
\end{equation}
Let $\overline{F}_k$ be a matrix induced by the action of $F_k$ on $\mathcal{V}^{N+1}/(\mathcal{K}+\mathcal{L})$, and we define the $q$-middle convolution $mc_{\lambda}$ as $\mathcal{E}\longrightarrow \mathcal{E}\ (E_{\boldsymbol{B},\boldsymbol{b}} \longmapsto E_{\boldsymbol{\overline{F}},\boldsymbol{b}})$.
\end{defi}

We abbreviated that modules $(\boldsymbol{B},\mathcal{V}),\, (\boldsymbol{F},\mathcal{V}^{N+1}),\, (\boldsymbol{\overline{F}},\mathcal{V}^{N+1}/(\mathcal{K}+\mathcal{L}))$ are $\mathcal{V},\, \mathcal{V}^{N+1},\, \mathcal{V}^{N+1}/(\mathcal{K}+\mathcal{L})$ respectively. Moreover, we set
\[ c_{\lambda}(\boldsymbol{B})=\boldsymbol{F},\ \ c_{\lambda}(\mathcal{V})=\mathcal{V}^{N+1},\ \ mc_{\lambda}(\boldsymbol{B})=\boldsymbol{\overline{F}},\ \ mc_{\lambda}(\mathcal{V})=\mathcal{V}^{N+1}/(\mathcal{K}+\mathcal{L}).\]

Here a $q$-analog of middle convolution was defined. We can also give an integral representation of $q$-convolution by $q$-analog of Euler transformation. We will describe it in detail in Section 2.

By the way, we would like to understand $q$-middle convolution as a transformation for the analog of Fuchsian equation. From now on, we set $q \in \mathbb{C},\, 0<|q|<1,\, \sigma _x:x \longmapsto qx$. We set a linear $q$-difference equation with polynomial coefficient
\begin{equation}\label{ea}
E_A\, :\ \sigma _xY(x)=A(x)Y(x),\ \ A(x)=\sum _{k=0}^NA_kx^k\ \ (A_k\in \mathrm{M}_m(\mathbb{C})).
\end{equation}
Moreover, we let $A_{\infty}=A_N.$ We define ``Fuchsian" $q$-difference equations.

\begin{defi}$\mathrm{(Fuchsian\ type\ equation)}$\ \ For an equation $E_A,$ if $A_0,\, A_{\infty}\in \mathrm{GL}_m(\mathbb{C})$, then we call $E_A$ Fuchsian type $q$-difference equation. 
\end{defi}

Although we cannot apply the $q$-middle convolution to this Fuchsian equation directly, we see that the equation $E_A$ is connected with $E_{\boldsymbol{B},\boldsymbol{b}}$ by simple transformations. We consider $m\times m$ matrix system $E_R$ with rational coefficients
\begin{equation}
E_R\, :\ \sigma _xY(x)=R(x)Y(x).
\end{equation}
As gauge transformations for the solution $Y(x)$ of the equation $E_R$, we consider only two types in this paper. The first one is the transformation
\begin{equation}
\varphi _P:Y(x)\longmapsto \widetilde{Y}(x)=PY(x)\ \ (P\in \mathrm{GL}_m(\mathbb{C})).
\end{equation}
The second one is the transformation
\begin{equation}
\varphi _f:Y(x) \longmapsto \widetilde{Y}(x)=f(x)Y(x),
\end{equation}
where $f(x)$ is solution of $\sigma _xf(x)=Q(x)f(x)\, $($Q(x)$ is a scalar rational function). This function $f(x)$ can be expressed by using the functions
\begin{equation*}
(ax;q)_{\infty},\ \ \ \ \ \ \vartheta _q (x).
\end{equation*}
Here we set
\begin{align*}
(a_1, \ldots ,a_n;q)_0&=1,\\
(a_1, \ldots ,a_n;q)_m&=\prod _{i=1}^n \prod _{j=0}^{m-1}(1-a_iq^j)\, (m \in \mathbb{Z}_{>0}),\\
(a_1, \ldots ,a_n;q)_{\infty }&=\lim_{m \rightarrow \infty}(a_1, \ldots ,a_n;q)_m,\\
\vartheta _q (x)&=\prod _{n=0}^{\infty}(1-q^{n+1})(1+xq^n)(1+x^{-1}q^{n+1}).
\end{align*}
To be specific, for the solution $Y(x)$ of the equation $E_R$,
\begin{align*}
&\text{if\ we\ put}\ \widetilde{Y}(x)=(ax;q)_{\infty}Y(x),\ \text{then}\ \sigma _x\widetilde{Y}(x)=(1-ax)R(x)\widetilde{Y}(x);\\
&\text{if\ we\ put}\ \widetilde{Y}(x)= \frac{1}{\vartheta _q (x)}Y(x),\ \text{then}\ \sigma _x\widetilde{Y}(x)=xR(x)\widetilde{Y}(x);\\
&\text{if\ we\ put}\ \widetilde{Y}(x)= \frac{\vartheta _q (x)}{\vartheta _q (ax)}Y(x)\, (a\in \mathbb{C}\backslash \{ 0 \} ),\ \text{then}\ \sigma _x\widetilde{Y}(x)=aR(x)\widetilde{Y}(x).
\end{align*}

We define a family of equations by modulo $\varphi _P$ and $\varphi _f$. We interpret the $q$-middle convolution as the transformation of the family of equations. From arbitrary equation $E_R$, we obtain $\Tilde{E}_R$:
\begin{align}
\Tilde{E}_R\, :\ \sigma _x\widetilde{Y}(x)&=A(x)\widetilde{Y}(x),\\
A(x)&=\sum_{i=0}^NA_ix^i\, (A_k\in \mathrm{M}(m,\mathbb{C}),\ A_0,A_N \ne 0,\ \forall a \in \mathbb{C}\, ;\, A(a)\ne 0),
\end{align}
which is determined up to multiplication of  constant and similarity transformations by $\varphi _P$.

We call $\Tilde{E}_R$ the canonical form of the equation $E_R$. In general case, for canonical form $\sigma _x\widetilde{Y}(x)=A(x)\widetilde{Y}(x)$ of $E_{\boldsymbol{B},\boldsymbol{b}}$, we obtain

\begin{align}
&A(x)=T(x)B(x),\ \ T(x)=\prod _{i=1}^N \Bigl( 1-\frac{x}{b_i}\Bigr) ,\\ 
&\label{rel1}A_0=1_m-B_0,\ A_{\infty}=b_{\infty}B_{\infty},\ B_0=1_m-\sum _{i=1}^NB_i-B_{\infty},\ b_{\infty}=\prod _{i=1}^N(-b_i^{-1}), \\
&\mathrm{rank}B_i = \begin{cases}
 m-n_1^k & (b_i=a_k\in Z_R)\\
 m & (b_i\notin Z_R)
\end{cases}\ (1 \le i \le N,\ n_1^k=\mathrm{dim\, ker}A(a_k)).
\end{align}

\begin{rem}
\normalfont{The }definition of the Fuchsian type equation may not be appropriate. We look at Heine's $q$-hypergeometric function
\begin{equation}
{}_2\varphi _1(\alpha ,\beta ,\gamma ;q;x)=\sum _{n=0}^{\infty}\frac{(\alpha ;q)_n(\beta ;q)_n}{(q;q)_n(\gamma ;q)_n}x^n.
\end{equation}
Here $u(x)={}_2\varphi _1(\alpha ,\beta ,\gamma ;q;x)$ satisfies the equation
\begin{equation}
\{ (1-\sigma _x)(1-q^{-1}\gamma \sigma _x)-x(1-\alpha \sigma _x)(1-\beta \sigma _x) \} u(x)=0.
\end{equation}
If we set $\displaystyle{v(x)=\frac{1}{x}\sigma _xu(x)}$ and $\displaystyle{Y(x)=\left( \begin{smallmatrix} u(x) \\ v(x) \end{smallmatrix} \right) ,}$ then we obtain
\begin{equation}\label{fu}
\sigma _xY(x)=\frac{1}{x(q\alpha \beta x-\gamma )}\begin{pmatrix} 0 & x^2(q\alpha \beta x-\gamma )\\ -x+1 & x\{ (\alpha +\beta )x-q^{-1}\gamma -1\} \end{pmatrix}Y(x).
\end{equation}
Although this is not Fuchsian $q$-difference equation in our sence, this equation transforms to Fuchsian type equation by a simple transformation:
\begin{equation}
Y(x)\longmapsto \Tilde{Y}(x)=\begin{pmatrix} 1 & 0\\ 1 & -x \end{pmatrix}Y(x)=\begin{pmatrix} u(x)\\ (1-\sigma _x)u(x)\end{pmatrix}.
\end{equation}
$\Tilde{Y}(x)$ satisfies Fuchsian $q$-difference equation
\begin{equation}\label{he}
\sigma _x\Tilde{Y}(x)=\frac{1}{\alpha \beta x-q^{-1}\gamma}\begin{pmatrix} \alpha \beta x-q^{-1}\gamma & -\alpha \beta x+q^{-1}\gamma \\ (1-\alpha )(1-\beta )x & (\alpha +\beta -\alpha \beta )x-1 \end{pmatrix}\Tilde{Y}(x).
\end{equation}
Although we do not introduce such transformations, Saloy used a transformation by rational component matrix as a gauge transformation in \cite{Sa}. We think that our Fuchsian $q$-difference equation corresponds to the Schlesinger normal form in the differential case. Although we do not call the equation (\ref{fu}) Fuchsian type, we might have to do. On the other hand, in the differential case, there exists Fuchsian differential equations which cannot be written in the Schlesinger normal form. We set $y_i(x)\, (i=1,2)$ the components of a solution $Y(x)$ of a equation
\begin{equation}\label{fu2}
\frac{dY}{dx}(x)=R(x)Y(x)\ \ \ \ (R(x)\ \mathrm{is\ rational\ function}).
\end{equation}
If singularities of $y_i(x)$ are at most regular singularities, we call the equation (\ref{fu2}) Fuchsian differential equation. Regular singularity is defined from local properties of solution. In more detail, if function $y(x)$ is not holomorphic at $x=x_0$ and for any $\underline{\theta},\overline{\theta}\, (\underline{\theta}<\overline{\theta}) ,$ there exists $n_0\in \mathbb{Z}_{>0}$ such that
\[ \lim _{\underline{\theta}<\arg (x-x_0)<\overline{\theta},\, x\rightarrow x_0}|x-x_0|^{n_0}|y(x)|=0,\]
we call $x=x_0$ the regular singularity of $y(x)$. Here we consider the equation of Schlesinger normal form
\[ \frac{dY}{dx}(x)=\left( \sum _{i=1}^N\frac{A_i}{x-a_i}\right) Y(x)\ \ \ \ (a_i\in \mathbb{C},\, A_i\in \mathrm{M}_m(\mathbb{C})),\]
that is, a special case of the Fuchsian differential equation.\ $\square$
\end{rem}

We can think that our Fuchsian type equation actually Fuchsian because Carmichael's theorem has been establish in \cite{Ca}.
\begin{theo}$\mathrm{(Carmichael)}$\ \ Let $\alpha _j^{\xi}\ (1\le j\le m,\, \xi =0,\infty )$ the eigenvalues of $A_{\xi}\in \mathrm{GL}_n(\mathbb{C})$, we assume further that $A_{\xi}$ are semi-simple and
\[ \frac{\alpha _j^{\xi}}{\alpha _k^{\xi}} \not\in q^{\mathbb{Z}_{>0}} \qquad (\forall j, \forall k).\] 
Then, there exist unique solutions $Y_{\xi}(x)$ of the equation (\ref{ea}) with the following properties,
\begin{eqnarray}
Y_0(x)&=&\widehat{Y}_0(x)x^{D_0},\\
Y_{\infty}(x)&=&q^{\frac{N}{2}u(u-1)}\widehat{Y}_{\infty}(x)x^{D_{\infty}},
\qquad
u=\frac{\log x}{\log q}.
\end{eqnarray}
Here $\widehat{Y}_{\xi}(x)$ is a holomorphic and invertible matrix at $x=\xi$ such that $\widehat{Y}_{\xi}(\xi )=C_{\xi}\in \mathrm{GL}(m,\mathbb{C})$ and $A_{\xi}=C_{\xi}q^{D_{\xi}}C_{\xi}^{-1},\, D_{\xi}=\mathrm{diag}(\log \alpha _j^{\xi}/\log q).$
\end{theo}

\begin{rem}
\normalfont{The }functions used in the above theorem
\begin{equation}
x^{\log \theta /\log q},\ \ \ \ q^{u(u-1)/2}\ \ (u=\log x/\log q)
\end{equation}
are solutions of the following equations, respectively,
\begin{equation}
\sigma _xf(x)=\theta f(x),\ \ \ \ \sigma _xf(x)=xf(x).
\end{equation}
Hence instead of these functions, we can use the following single-valued functions as solutions of the above equations,
\begin{equation}
\frac{\vartheta _q (x)}{\vartheta _q (\theta x)},\ \ \frac{1}{\vartheta _q (x)}.
\end{equation}
These functions are widely used in recent years, we use these in this paper.\ $\square$
\end{rem}

The purpose of this study is to examine properties of the $q$-middle convolution. Let us describe the contents of this paper. In the 2nd section, we show that $q$-convolution can be expressed by a $q$-analog of Euler transformation. In the 3rd section, we define the spectral type and the rigidity index for the equation $E_R$. Spectral types are defined by the size of Jordan cells of $A_0,\, A_{\infty}$ and types of elementary divisors of $A(x)$.  Notice that the rigidity index is not only determined by data of $B_k$ of $B(x)$, but also by data of elementary divisors of coefficient polynomial $A(x)$ of canonical form $\Tilde{E}_R$. In the 4th section, we prove the three main theorems.

\begin{theo}\label{Fuchs}$(\mathrm{Fuchsian\ type\ equation})$\ \ If equation $E_R$ is Fuchsian type equation, then $mc_{\lambda}(E_R)$ is also Fuchsian type equation. 
\end{theo}
Here we assume that next conditions $(\ast ),(\ast \ast)$ after the manner of Dettweiler and Reiter in \cite{DR1}. (These conditions are generally satisfied if $\mathrm{dim}\mathcal{V}=1$ or $\mathrm{dim}\mathcal{V}>1$ and $\boldsymbol{B}$ is irreducible)

\begin{defi}\label{ast}\ \ We define the conditions $(\ast ),(\ast \ast)${\rm :}
\begin{align*}
(\ast )\, :\ &\forall i \in \{ 0, \ldots ,N \} , \forall \tau \in \mathbb{C}\ ;\ \  \bigcap _{j \neq i}\ker B_j\! \cap \! \ker (B_i\! +\! \tau 1_m)\! =\! 0,\\
(\ast \ast )\, :\ &\forall i \in \{ 0, \ldots ,N \} , \forall \tau \in \mathbb{C}\ ;\ \  \sum _{j \neq i}\mathrm{im}B_j\! +\! \mathrm{im}(B_i\! +\! \tau 1_m)\! =\! \mathcal{V}.
\end{align*}
\end{defi}
\begin{theo}\label{irre}$(\mathrm{irreducibility})$\ \ If $(\ast ),(\ast \ast )$ are satisfied, then $\mathcal{V}$ is irreducible if and only if $mc_{\lambda}(\mathcal{V})$ is irreducible.
\end{theo}
\begin{theo}\label{index}$(\mathrm{rigidity\ index})$\ \ If $(\ast ) ,(\ast \ast )$ are satisfied, then $mc_{\lambda}$ preserves rigidity index of Fuchsian equation $E_R$.
\end{theo}

To prove these theorems, we do not need for the following conditions in the Theorem 1.6 :
\[ A_0, A_{\infty}\, :\, \mbox{semi-simple},\ \ \ \ \frac{\theta _j}{\theta _k}, \frac{\kappa _j}{\kappa _k}
 \not\in q^{\mathbb{Z}_{>0}}\ (\theta _i, \kappa _i:\, \mathrm{eigenvalues\ of}\ A_0,\, A_{\infty}\, \mathrm{respectively}).\]
We will explain ``rigidity indexh in the section 3. This is defined by ``spestral type'' of the Fuchsian equation $E_R$. 

\section{Integral representation of $q$-convolution}

We gave a $q$-analog of convolution as a transformation of the $q$-difference equations. We can also give an integral representation of ``$q$-convolution'' by a $q$-analog of Euler transformation. In this section, we show

\begin{theo}\ \ For\ the\ solution\ $Y(x)$\ of\ the\ equation\ $E_{\boldsymbol{B},\boldsymbol{b}}$,\ let\ $\widehat{Y}(x)={}^t({}^t\widehat{Y}_0(x), \ldots ,{}^t\widehat{Y}_N(x))$\ by
\begin{equation}
\widehat{Y}_i(x)=\int _0^{\infty}\frac{P_{\lambda}(x,s)}{s-b_i}Y(s)\ d_qs,\ b_0\! =0,\ P_{\lambda}(x,s)=\frac{(q^{\lambda +1}sx^{-1};q)_{\infty}}{(qsx^{-1};q)_{\infty}}.
\end{equation}
Then,\ $\widehat{Y}(x)$\ is\ the\ solution\ of\ the\ equation\ $E_{\boldsymbol{F},\boldsymbol{b}}${\rm (}see Definition \ref{conv}{\rm )}. Here Jackson integral is defined by
\[ \int _0^{\infty}f(x)\, d_qx=(1-q)\sum _{n=-\infty}^{\infty}q^nf(q^n).\]
\end{theo}

$Proof.$\ \ $P_{\lambda}(x,s)$ \rm{is a solution of partial difference equations}
\begin{equation*}
(\sigma _x-\sigma _s^{-1})y(x,s)=0,\ \sigma _xy(x,s)=\frac{1-q^{\lambda}sx^{-1}}{1-sx^{-1}}y(x,s).
\end{equation*}
Hence $P_{\lambda}(x,s)$ is a solution of 
\begin{equation*}
\frac{\sigma _xP_{\lambda}(x,s)}{s-b_i}=\frac{x-q^{\lambda}b_i}{x-b_i}\frac{P_{\lambda}(x,s)}{s-b_i}+\frac{x}{x-b_i}\frac{\sigma _s^{-1}-1}{s}P_{\lambda}(x,s).
\end{equation*}
Moreover, this function is independent to $b_i\in \mathbb{C}$. By multiplying $Y(s)$, and by Jackson integral calculation, we obtain
\begin{equation*}
\sigma _x\widehat{Y}_i(x)=\left\{ 1+\frac{(1-q^{\lambda})b_i}{x-b_i} \right\} \widehat{Y}_i(x)+\frac{x}{x-b_i}\int _0^{\infty}\frac{\sigma _s^{-1}-1}{s}P_{\lambda}(x,s)Y(s)\ d_qs.
\end{equation*}
Meanwhile, we obtain
\begin{align*}
\int _0^{\infty}& \frac{\sigma _s^{-1}-1}{s}P_{\lambda}(x,s)\cdot Y(s)\ d_qs\\
&=\int _0^{\infty}P_{\lambda}(x,s)\frac{1}{s}\{ \sigma _s Y(s)-Y(s) \} \ d_qs\\
&=\int _0^{\infty}P_{\lambda}(x,s)\frac{1}{s}\biggl( B_{\infty}+\sum _{j=1}^N \frac{B_j}{1-\frac{s}{b_j}}-1_m \biggr) Y(s)\ d_qs\\
&=\int _0^{\infty}P_{\lambda}(x,s)\biggl\{ \frac{1}{s}\biggl( B_{\infty}+\sum _{j=1}^NB_j-1_m\biggr) -\sum _{j=1}^N\frac{1}{s-b_j}B_j\biggr\} Y(s)\ d_qs\\
&=-\int _0^{\infty}P_{\lambda}(x,s)\sum _{j=0}^N\frac{1}{s-b_j}B_j\cdot Y(s)\ d_qs\ \biggl( b_0=0,\ B_0=1_m-\sum _{i=1}^NB_i-B_{\infty}\biggr) \\
&=-\sum _{j=0}^NB_j\int _0^{\infty}\frac{P_{\lambda}(x,s)}{s-b_j}Y(s)\ d_qs\\
&=-\sum _{j=0}^NB_j\widehat{Y}_j(x).
\end{align*}
Here $\widehat{Y}_i(x)$ satisfies
\begin{align*}
\sigma _x\widehat{Y}_i(x)&=\left\{ 1+\frac{(1-q^{\lambda})b_i}{x-b_i} \right\} \widehat{Y}_i(x)-\frac{x}{x-b_i}\sum _{j=0}^NB_j\widehat{Y}_j(x)\\
&=\widehat{Y}_i(x)-\sum _{j=0}^NB_j\widehat{Y}_j(x)+\frac{1}{1-\frac{x}{b_i}}\biggl\{ -(1-q^{\lambda})\widehat{Y}_i(x)+\sum _{j=0}^NB_j\widehat{Y}_j(x) \biggr\} .
\end{align*}
Therefore, $\widehat{Y}(x)$ is a solution of the equation $E_{\boldsymbol{F},\boldsymbol{b}}$.\ $\square$

From the above, we proved that $q$-convolution can be expressed by a $q$-analog of Euler transformation.

\section{Rigidity index of $q$-difference equations}

In this section, we define the spectral type and the rigidity index of the equation $E_R$. We set the coefficient $A(x)=\sum _{k=0}^NA_kx^k$ of the canonical form of a Fuchsian equation $E_R$.

\begin{defi}
Let\ $A_{\xi} \sim \bigoplus _{i=1}^{l_{\xi}}\bigoplus _{j=1}^{s_i^{\xi}} J(\alpha _i^{\xi},t_{i,\, j}^{\xi})\ (J(\alpha ,t):\, Jordan\ cell,\ t_{i,\, j+1}^{\xi} \le t_{i,\, j}^{\xi}).$ Moreover,\ let\ $\{ m_{i,\, k}^{\xi}\} _k$\ denote\ the\ conjugate\ of\ $\{ t_{i,\, j}^{\xi}\} _j$\ in\ Young\ diagram. We\ call
\[ S_{\xi}:m_{1,1}^{\xi} \ldots m_{1,{t_{1,1}^{\xi}}}^{\xi}, \ldots ,m_{{l_{\xi}},1}^{\xi} \ldots m_{{l_{\xi}},{t_{{l_{\xi}},1}^{\xi}}}^{\xi}\]
spectral\ type\ of\ $A_{\xi}.$
\end{defi}

\begin{defi}
Let\ $Z_{A}=\{a \in \mathbb{C}\, ;\, \mathrm{det}\, A(a)=0 \}$\ and\ denote\ by\ $d_i\, (1\le i\le m)$ the elementary\ divisors\ of\ $\mathrm{det}\, A(x)\, (d_{i+1}|d_i)$. For\ any\ $a_i\in Z_A$,\ we\ denote\ by\ $\{ \Tilde{n}_k^i\} _k$\ the\ orders\ of\ zeros\ $a_i$\ of\ $\{ d_k\} _k$. We set $\{ n_j^i\} _j$\ the\ conjugate\ of\ $\{ \Tilde{n}_k^i\} _k$.\ We\ call
\[ S_{\mathrm{div}}:n_1^1 \ldots n_{k_1}^1, \ldots ,n_1^l \ldots n_{k_l}^l\]
spectral\ type\ of\ $A(x)$.
\end{defi}

\begin{defi}
We\ call\ $S(E_R)=(S_0;S_{\infty};S_{\mathrm{div}})$\ spectral\ type\ of\ $E_R.$
\end{defi}

From the above, we define the rigidity index.

\begin{defi}
We\ define\ the\ rigidity\ index\ $\mathrm{idx}(E_R)$\ of\ the\ equation\ $E_R$\ as
\begin{equation}
\mathrm{idx}(E_R)=\ \sum_{\xi =0,\infty}\sum_{i=1}^{l_{\xi}}\sum_{j=1}^{t_{i,1}^{\xi}}(m_{i,\, j}^{\xi})^2+\sum _{i=1}^l\sum_{j=1}^{k_i}(n_j^i)^2-m^2N.
\end{equation}
\end{defi}

For example, we consider
\begin{align*}
&E_1\, :\ \sigma _xY(x)=A(x)Y(x),\ \ \ A(x)=A_0+A_1x+A_{\infty}x^2,\\
&A_0\sim J(\alpha _1^0,2)\oplus J(\alpha _1^0,1)^{\oplus 2}\oplus J(\alpha _2^0,1),\ A_{\infty}\sim J(\alpha _1^{\infty},1)^{\oplus 3}\oplus J(\alpha _2^{\infty},1)\oplus J(\alpha _3^{\infty},1),\\
&A(x)\sim \mathrm{diag}((x-a_1)(x-a_2)^2(x-a_3)(x-a_4),\, (x-a_1)(x-a_2),\, (x-a_1)(x-a_2),\, x-a_1,1)\\
&\ (\alpha _i^0\ne \alpha _j^0,\, \alpha _i^{\infty}\ne \alpha _j^{\infty},\, a_i\ne a_j\, (i\ne j)).
\end{align*}
Spectral type and rigidity index of the equation $E_1$ are
\[ S(E_1)\, :\ 31,\! 1;3,\! 1,\! 1;4,\! 31,\! 1,\! 1,\ \ \ \ \mathrm{idx}(E_1)=0.\]

\begin{rem}
\normalfont{We}\ can\ also\ express\ the\ rigidity\ index\ $\mathrm{idx}(E_R)$\ of\ the\ equation\ $E_R$\ as
\begin{equation}
\mathrm{idx}(E_R)=\mathrm{dim}Z(A_0)+\mathrm{dim}Z(A_{\infty})+\sum _{i=1}^l\sum_{j=1}^{k_i}(n_j^i)^2-m^2N.
\end{equation}
Here, we let $Z(A)=\{ X\in \mathrm{GL}_m(\mathbb{C})\, ;\, AX=XA\} \, (A\in \mathrm{M}_m(\mathbb{C}))$.\ $\square$
\end{rem}

We can easily check the next facts.

\begin{prop}
\begin{align*}
\mathrm{(i)}\ \ &\sum_{i=1}^{l_{\xi}}\sum_{j=1}^{t_{i,1}^{\xi}}m_{i,\, j}^{\xi}=m,\ \ \sum _{i=1}^l\sum_{j=1}^{k_i}n_j^i=Nm.\\
\mathrm{(ii)}\ \ &n_i=\textstyle{\sum_{j=1}^{k_i}}n_j^i\, is\ a\ multiplicity\ of\ \mathrm{det}A(x)\ of\ zeros\ a_i\in Z_A.\\
\mathrm{(iii)}\ \ &\mathop{\rm idx} (E_R)\ is\ even\ number.
\end{align*}
\end{prop}

After the definition of $q$-analog of spectral type and rigidity index, let's look at some examples. At first, we consider the Heine's $q$-hypergeometric equation $E_2\,$:\ (\ref{he}). It is easy to confirm that the equation $E_2$ has generally the following data:

\begin{equation}
S(E_2)\, :\ 1,\! 1;1,\! 1;1,\! 1,\ \ \ \ \mathrm{idx}(E_2)=2\rm{.}
\end{equation}

Moreover, we consider generalized $q$-hypergeometric equation
\begin{align}
&E_3\, :\ \sigma _xY(x)=A(x)Y(x),\ \ A(x)=\begin{pmatrix} 0 & f_0& & \\
 & \ddots & \ddots & \\
 & & 0 & f_0\\
 -f_m & \cdots &-f_2 & -f_1 \end{pmatrix},\\
&\ f_0\sigma _x^m+f_1\sigma _x^{m-1}+\cdots +f_m=\prod _{k=1}^m\left( \frac{b_k}{q}\sigma _x-1\right) -\lambda x\prod _{k=1}^m(a_k\sigma _x-1)\ (a_k,b_k,\lambda \in \mathbb{C}^{\ast}).
\end{align}
We set $A(x)=A_0+A_{\infty}x\ (A_k\in \mathrm{M}_m(\mathbb{C}))$. We obtain the data of the equation $E_3$ as
\begin{align}
&\mathrm{Ev}(A_0)=\left\{ \frac{q}{b_1},\, \ldots ,\, \frac{q}{b_m}\right\} ,\ \ \mathrm{Ev}(A_{\infty})=\left\{ \frac{1}{a_1},\, \ldots ,\, \frac{1}{a_m}\right\} ,\\
&\mathrm{zeros\ of\ det}A(x)\ \mathrm{are}\ \frac{1}{\lambda}\ \mathrm{and}\ \frac{1}{\lambda}\prod _{k=1}^m\frac{b_k}{qa_k}\ (\mathrm{multiplicity}: m-1).
\end{align}
Here we denote by $\mathrm{Ev}(A_{\xi})\, (\xi =0,\infty )$ the set of eigenvalues of $A_{\xi}$. Therefore, we generally obtain rigidity index of the equation $E_3$ as
\[ \mathrm{idx}(E_3)=1^2\times m+1^2\times m+1^2+(m-1)^2-1\times m^2=2.\]

\begin{rem}
\normalfont{In general case, since} we can also express the Fuchsian equation $E_A\, :\, \sigma _x^{-1}Y(x)=A(q^{-1}x)^{-1}Y(x),$ we expect $\mathrm{idx}(E_{A^{-1}})=\mathrm{idx}(E_A)$. Let us check this fact. We put
\[ \Tilde{A}(x)=\mathrm{det}A(x)\, A(x)=\sum _{k=0}^{N(m-1)}\Tilde{A}_kx^k,\ \ \ \Tilde{A}_{\infty}=\Tilde{A}_{N(m-1)},\]
then we get 
\[ A_0\Tilde{A}_0=1_m,\ \ \ \ A_{\infty}\Tilde{A}_{\infty}=\kappa 1_m\ (\kappa \in \mathbb{C}\backslash \{ 0\} ).\]
Moreover, the spectral type $S(E_{A^{-1}})=(\Tilde{S}_0;\Tilde{S}_{\infty};\Tilde{S}_{\mathrm{div}})$ satisfies $\Tilde{S}_0=S_0,\, \Tilde{S}_{\infty}=S_{\infty}$ and 
\[ \Tilde{S}_{\mathrm{div}}\, :\, \underbrace{m\ldots m}_{n_1-k_1}m-n^1_{k_1}\ldots m-n^1_1,\, \ldots ,\, \underbrace{m\ldots m}_{n_l-k_l}m-n^l_{k_l}\ldots m-n^l_1\] 
because $\Tilde{A}(x)\sim \mathrm{det}A(x)\, \mathrm{diag}(d_i^{-1}).$ Therefore, we obtain
\begin{align*}
\mathrm{idx}(E_{A^{-1}})&=\mathrm{dim}Z(\Tilde{A}_0)+\mathrm{dim}Z(\Tilde{A}_{\infty})+\sum _{i=1}^l\left\{ m^2(n_i-k_i)+\sum _{j=1}^{k_i}(m-n_j^i)^2\right\} -N(m-1)m^2\\
&=\mathrm{dim}Z(A_0)+\mathrm{dim}Z(A_{\infty})+\sum _{i=1}^l\left\{ m^2n_i-2m\sum _{j=1}^{k_i}n_j^i+\sum _{j=1}^{k_i}(n_j^i)^2\right\} -N(m-1)m^2\\
&=\mathrm{dim}Z(A_0)+\mathrm{dim}Z(A_{\infty})+(m^2-2m)\cdot Nm+\sum _{i=1}^l\sum _{j=1}^{k_i}(n_j^i)^2-N(m-1)m^2\\
&=\mathrm{dim}Z(A_0)+\mathrm{dim}Z(A_{\infty})+\sum _{i=1}^l\sum _{j=1}^{k_i}(n_j^i)^2-Nm^2\\
&=\mathrm{idx}(E_A).\ \square
\end{align*}
\end{rem}

In the next section, we study how these data are changed by $q$-middle convolution in detail.

\section{Properties of $q$-middle convolution}

In this section, we prove the three theorems. \\
\\
\textbf{Theorem \ref{Fuchs}}\ (Fuchsian type equation)\ \ \it{If equation $E_R$ is Fuchsian type equation, then $mc_{\lambda}(E_R)$ is also Fuchsian type equation}\rm{.} \\
\\
\rm{\textbf{Theorem \ref{irre}}\ (irreducibility)}\ \ \it{If $(\ast ),(\ast \ast )$ are satisfied, then $\mathcal{V}$ is irreducible if and only if $mc_{\lambda}(\mathcal{V})$ is irreducible}\rm{.}\\
\\
\rm{\textbf{Theorem \ref{index}}\ (rigidity index)}\ \ \it{If $(\ast ) ,(\ast \ast )$ are satisfied, then $mc_{\lambda}$ preserves rigidity index of Fuchsian equation $E_R$}\rm{.}\\

About $(\ast ),(\ast \ast ),$ see Definition \ref{ast}. Theorem \ref{Fuchs} is proved easily by examining coefficient polynomial of canonical form of $c_{\lambda}(\Tilde{E}_R)$. Although many preparations are necessary for us to prove Theorem \ref{irre}, the outline is the same as method of Detteweiler and Reiter in \cite{DR1}. Finally, Theorem \ref{index} is proved by investigating in detail the change of spectral type of the equation $E_R$.

\subsection{Proof of Theorem \ref{Fuchs}.}

Here we prove the next theorem.\\
\\
\textbf{Theorem \ref{Fuchs}}\ (Fuchsian type equation)\ \ \it{If equation $E_R$ is Fuchsian type equation, then $mc_{\lambda}(E_R)$ is also Fuchsian type equation}\rm{.}\\

$Proof.$\ \ We put coefficients $A(x)=\sum _{k=0}^NA_kx^k\, (A_0,A_{\infty}\in \mathrm{GL}_m(\mathbb{C})),\ G(x)=\sum _{k=0}^NG_kx^k$ of canonical form of $E_{\boldsymbol{B},\boldsymbol{b}},\, E_{\boldsymbol{F},\boldsymbol{b}}\, (\boldsymbol{F}=c_{\lambda}(\boldsymbol{B}))$. From the relations (\ref{rel1}):
\begin{equation*}
A_0=1_m-B_0,\ A_{\infty}=b_{\infty}B_{\infty},\ B_0=1_m-\sum _{i=1}^NB_i-B_{\infty},\ b_{\infty}=\prod _{i=1}^N(-b_i^{-1})\ne 0,
\end{equation*}
we obtain $B_0-1_m,B_{\infty}\in \mathrm{GL}_m(\mathbb{C})$. For any 
$v={}^t({}^tv_0, \ldots ,{}^tv_N)\in \ker F_{\infty}\, (v_k\in \mathcal{V}),$ we get $G_{\infty}\in \mathrm{GL}_{(N+1)m}(\mathbb{C})$ because
\[ 0=G_{\infty}v=b_{\infty}F_{\infty}v=b_{\infty}{}^t({}^t(B_{\infty}s), \ldots ,{}^t(B_{\infty}s))\, (s=\textstyle{\sum _{i=0}^N} B_iv_i).\]

Meanwhile, for any $v={}^t({}^tv_0, \ldots ,{}^tv_N) \in \ker G_0, $ since
\[ 0=G_0v=(1_{(N+1)m}-F_0)v=(\textstyle{\sum _{i=1}^N}F_i+F_{\infty})v=(\textstyle{\sum _{i=1}^N}F_i+1_{(N+1)m}-\widehat{F})v,\]
we obtain $v=0.$ Hence $G_0\in \mathrm{GL}_{(N+1)m}(\mathbb{C})$. Therefore, $mc_{\lambda}(E_R)$ is a Fuchsian type equation.\ $\square$\\

\subsection{Proof of Theorem \ref{irre}.}

Here we derive a dimension formula of $q$-middle convolution. Moreover, we prove that $q$-middle convolution preserves irreducibility of the equation. The outline is the same as calculations of Detteweiler and Reiter in \cite{DR1}.

At first, linear spaces $\mathcal{K},\mathcal{L}$ satisfy the next proposition.

\begin{prop}\ \ $\mathcal{K},\mathcal{L}$ are $\boldsymbol{F}$-invariant subspaces of $\mathcal{V}^{N+1}.$
\end{prop}

$Proof.$\ \ (i)\ Let $J=\{ 1, \ldots ,N \}.$ For any $v={}^t({}^tv_0, \ldots ,{}^tv_N) \in \mathcal{K}\, (v_k\in \mathrm{ker}B_k)$, we get
\[ F_jv={}^t(0, \ldots ,\stackrel{j+1}{\stackrel{\vee}{(q^{\lambda}-1){}^tv_j}} , \ldots , 0)\in \mathcal{K}\ \ (j \in J).\]
Hence $F_j\mathcal{K}$ is subspace of $\mathcal{K}$. In the meantime, $F_{\infty}\mathcal{K}$ is subspace of $\mathcal{K}$ because for any $v \in \mathcal{K}$, we obtain $F_{\infty}v=(1_{(N+1)m}-\widehat{F})v=v \in \mathcal{K}$. Therefore, $\mathcal{K}$ is $\boldsymbol{F}$-invariant subspace of $\mathcal{V}^{N+1}.$

(ii)\ Let
\[ 1_{m,k}=\{ \delta _{i,k+1}\delta _{j,k+1}1_m\} _{1\le i,j\le N+1}=\mathrm{diag}(0, \ldots ,\stackrel{k+1}{\stackrel{\vee}{1_m}} , \ldots ,0).\]
For any $v \in \mathcal{L}$, we get 
\[ (\widehat{F}-(1-q^{\lambda})1_{(N+1)m})F_jv=(\widehat{F}-(1-q^{\lambda})1_{(N+1)m})1_{m,j}(\widehat{F}-(1-q^{\lambda})1_{(N+1)m})v=0\ \ (j \in J).\]
Hence $F_j\mathcal{L}$ is subspace of $\mathcal{L}$. Moreover, $F_{\infty}\mathcal{L}$ is subspace of $\mathcal{L}$ because for any $v \in \mathcal{L}$, we obtain 
\[ (\widehat{F}-(1-q^{\lambda})1_{(N+1)m})F_{\infty}v=(\widehat{F}-(1-q^{\lambda})1_{(N+1)m})(1_{(N+1)m}-\widehat{F})v=0.\]
Therefore, $\mathcal{L}$ is $\boldsymbol{F}$-invariant subspace of $\mathcal{V}^{N+1}.$\ $\square$
\\

The next facts are important as ``dimension formula".

\begin{prop}{\ }\\
$\mathrm{(i)}$\ If $\lambda = 0,$ then $\mathcal{K}$ is subspace of $\mathcal{L}$ and satisfies
\begin{equation} 
\mathcal{L}= \{ {}^t({}^tv_0, \ldots ,{}^tv_N) ; \textstyle{\sum _{j=0}^N} B_jv_j = 0 \} .
\end{equation}
$\mathrm{(ii)}$\ If $\lambda \neq 0,$ then $\mathcal{K}\cap \mathcal{L}=0, \mathcal{L}= \{ {}^t({}^th, \ldots ,{}^th) ;h \! \in \ker (A_{\infty}-q^{\lambda}b_{\infty}1_m) \}$ and
\begin{equation}
\dim (mc_{\lambda}(\mathcal{V}))=(N+1)m-\sum _{i=1}^N \dim \ker B_i-\dim \ker (A_0-1_m)-\dim \ker (A_{\infty}-q^{\lambda}b_{\infty}1_m).
\end{equation}
\end{prop}

$Proof.$\ \ (i)\ If $\lambda =0$, then $\mathcal{L}=\ker \widehat{F}.$ Here for any $v \in \mathcal{K},$ we obtain $\widehat{F}v=0$. Hence $v \in \mathcal{L}.$ Moreover, we obtain $\mathcal{L}= \{ {}^t({}^tv_0, \ldots ,{}^tv_N) ; \sum _{j=0}^N B_jv_j = 0 \} .$

(ii)\ If $\lambda \ne 0,$ for any $v \in \mathcal{K}\cap \mathcal{L},$ we obtain 
\[ 0=(\widehat{F}-(1-q^{\lambda})1_{(N+1)m})v=\widehat{F}v-(1-q^{\lambda})v=(q^{\lambda}-1)v.\] 
Hence we get $v=0$. For any $v={}^t({}^tv_0, \ldots ,{}^tv_N) \in \mathcal{L},$ we obtain $\widehat{F}v=(1-q^{\lambda})v.$ Consequently, we see $\sum _{j=0}^NB_jv_j=(1-q^{\lambda})v_i\, (i \in I=\{ 0, \ldots ,N \} ).$ Here $v_0=\dotsi =v_N$ and 
\[ \mathcal{L}= \{ {}^t({}^th, \ldots ,{}^th) ;h \! \in \ker (A_{\infty}-q^{\lambda}b_{\infty}1_m) \} .\] 
Therefore, we can compute $\dim (mc_{\lambda}(\mathcal{V}))\, $:
\begin{align*}
\dim (mc_{\lambda}(\mathcal{V}))&=\dim (\mathcal{V}^{N+1}/(\mathcal{K}+\mathcal{L}))\\
&=\dim (\mathcal{V}^{N+1})-\dim (\mathcal{K}+\mathcal{L})\\
&=\dim (\mathcal{V}^{N+1})-\dim \mathcal{K}-\dim \mathcal{L}\ \ (\because \ \mathcal{K}\cap \mathcal{L}=0)\\
&=(N+1)m-\sum _{i=0}^N\dim \ker B_i-\dim \ker (B_{\infty}-q^{\lambda}1_m)\\
&=(N+1)m-\sum _{i=1}^N\dim \ker B_i-\dim \ker (A_0-1_m)-\dim \ker (A_{\infty}-q^{\lambda}b_{\infty}1_m).\ \ \square 
\end{align*}

\begin{prop}\ \ If $\mathcal{W}$ is $\boldsymbol{B}$-invariant subspace of $\mathcal{V}$, then $\mathcal{W}^{N+1}$ is $\boldsymbol{F}$-invariant subspace. Moreover, $mc_{\lambda}(\mathcal{W})$ is submodule of $mc_{\lambda}(\mathcal{V}).$
\end{prop}

$Proof.$\ \ For any $w={}^t({}^tw_0, \ldots ,{}^tw_N) \in \mathcal{W}^{N+1}$ and $j \in J=\{ 1, \ldots ,N \}$, it is clear that
\[ F_jw={}^t(0, \ldots ,\stackrel{j+1}{\stackrel{\vee}{\textstyle{\sum _{i=0}^N}{}^t(B_iw_i)-(1-q^{\lambda}){}^tw_j}},\ldots ,0) \in \mathcal{W}^{N+1}.\]
Since $F_{\infty}w=(1_{(N+1)m}-\widehat{F})w=w-\widehat{F}w \in \mathcal{W}^{N+1},\ \mathcal{W}^{N+1}$ is $\boldsymbol{F}$-invariant subspace of $\mathcal{V}^{N+1}$. The second claim follows from 
\begin{equation}\label{rel2}
\mathcal{W}^{N+1} \cap (\mathcal{K}_{\mathcal{V}}+\mathcal{L}_{\mathcal{V}})=\mathcal{K}_{\mathcal{W}}+\mathcal{L}_{\mathcal{W}}.
\end{equation}
Hence we prove (\ref{rel2}). If $\lambda =0,\ \mathcal{K}$ is subspace of $\mathcal{L}$. If $\lambda \ne 0,$ then
\[ \mathcal{K}_{\mathcal{W}}+\mathcal{L}_{\mathcal{W}}=\mathcal{K}_{\mathcal{V}\cap \mathcal{W}}+\mathcal{L}_{\mathcal{V}\cap \mathcal{W}}\]
is subspace of $\mathcal{W}^{N+1} \cap (\mathcal{K}_{\mathcal{V}}+\mathcal{L}_{\mathcal{V}})$. Moreover, for any $w={}^t({}^tw_0, \ldots ,{}^tw_N) \in \mathcal{W}^{N+1} \cap (\mathcal{K}_{\mathcal{V}}+\mathcal{L}_{\mathcal{V}})$ and $i \in I=\{ 0, \ldots ,N \}$, we can let
\[ w_i=k_i+h\, (k_i \in \ker B_i,h \in \ker (A_{\infty}-q^{\lambda}b_{\infty}1_m)).\]
Here we obtain $\mathcal{W}\ni \sum _{i=0}^NB_iw_i= \sum _{i=0}^NB_i(k_i+h)=(1-q^{\lambda})h.$ Consequently, $h \in \mathcal{W}$. Moreover, we find $w \in \mathcal{K}_{\mathcal{W}}+\mathcal{L}_{\mathcal{W}}$ from $k_i=w_i-h \in \mathcal{W}$. Therefore, $\mathcal{W}^{N+1} \cap (\mathcal{K}_{\mathcal{V}}+\mathcal{L}_{\mathcal{V}})$ is subspace of $\mathcal{K}_{\mathcal{W}}+\mathcal{L}_{\mathcal{W}}.\ \square$\\

From now on, we assume the conditions $(\ast ),(\ast \ast)$. Here we can prove

\begin{prop}\label{404}\ \ If $(\ast \ast)$ is satisfied, then $mc_0(\mathcal{V}) \simeq \mathcal{V}.$
\end{prop}

$Proof.$\ \ If $\lambda =0,$ then we get $\mathcal{K}+\mathcal{L}= \mathcal{L}= \{ {}^t({}^tv_0, \ldots ,{}^tv_N) ; \sum _{j=0}^NB_jv_j=0\}$. Let 
\[ \phi :{}^t({}^tv_0, \ldots ,{}^tv_N) \longmapsto \sum _{j=0}^NB_jv_j.\]
Then $\phi :\mathcal{V}^m \longrightarrow \mathcal{V}$ is surjection from a condition $(\ast \ast)$. For any $v={}^t({}^tv_0, \ldots ,{}^tv_N) \in \mathcal{V}^{N+1}$, we get
\begin{align*}
(\phi \circ F_j)(v)&=\phi ({}^t(0, \ldots ,\stackrel{j+1}{\stackrel{\vee}{{}^ts}} , \ldots , 0))=B_js=(B_j \circ \phi )(v),\ s=\textstyle{\sum _{i=0}^N}B_iv_i\, (j \in J=\{ 1, \ldots ,N \} ),\\
(\phi \circ F_{\infty})(v)&=(\phi \circ (1_{(N+1)m}-\widehat{F}))(v)=\phi ({}^t({}^tv_0-{}^ts, \ldots ,{}^tv_N-{}^ts))=B_{\infty}s=(B_{\infty} \circ \phi )(v).
\end{align*}
Therefore, we obtain 
\begin{equation*}
\mathcal{V}= \mathrm{im}(\phi ) \simeq \mathcal{V}^{N+1}/{\ker (\phi )} = \mathcal{V}^{N+1}/(\mathcal{K}+\mathcal{L}) = mc_0(\mathcal{V}).\ \square 
\end{equation*}

Here we introduce a transformation $\psi _{\mu}$ in expedient.

\begin{defi}\ \ For $\boldsymbol{T}=(T_1, \ldots ,T_N,T_{\infty}) \in (\mathrm{M}_{(N+1)m}(\mathbb{C}))^{N+1},$ we define
\begin{equation}
\psi _{\mu}:(\mathrm{M}_{(N+1)m}(\mathbb{C}))^{N+1} \longrightarrow (\mathrm{M}_{(N+1)m}(\mathbb{C}))^{N+1},\ (T_1, \ldots ,T_N,T_{\infty}) \longmapsto (T_1, \ldots ,T_N,T_{\infty}+\mu 1_{(N+1)m}).
\end{equation}
We set the module $\psi _{\mu}(\mathcal{V})=(\psi _{\mu}(\boldsymbol{T}),\mathcal{V}).$
\end{defi}

Here $\psi _{\mu}$ preserves irreducibility of equations clearly. Moreover, we introduce a transformation $\Psi _{\lambda}$.

\begin{defi}\ \ We define $\Psi _{\lambda}:\mathcal{E}\longrightarrow \mathcal{E}$,
\begin{equation}
\Psi _{\lambda}=\psi _{1-q^{\lambda}} \circ c_{\lambda}.
\end{equation}
Let $\widetilde{\boldsymbol{F}}=\Psi _{\lambda}(\boldsymbol{B}),\ \Psi _{\lambda}(\mathcal{V})=(\widetilde{\boldsymbol{F}},\mathcal{V}^{N+1}).$ We let $\check{F}_k$ be a matrix induced by the action of $F_k$ on $\mathcal{V}^{N+1}/(\mathcal{K}+\mathcal{L}).$ Moreover, we define $\overline{\Psi}_{\lambda}:\mathcal{E}\longrightarrow \mathcal{E}$,
\begin{equation}
\overline{\Psi}_{\lambda}(\boldsymbol{B})=\check{\boldsymbol{F}},\ \overline{\Psi}_{\lambda}(\mathcal{V})=\mathcal{V}^{N+1}/(\mathcal{K}+\mathcal{L})=(\check{\boldsymbol{F}},\mathcal{V}^{N+1}/(\mathcal{K}+\mathcal{L})).
\end{equation}
\end{defi}

Here the following facts are proved in the same way as above.

\begin{prop}\ \ $\mathcal{K},\mathcal{L}$ are $\widetilde{\boldsymbol{F}}$-invariant.
\end{prop}

\begin{prop}\ \ If $\mathcal{W}$ is $\boldsymbol{B}$-invariant subspace of $\mathcal{V}$, then $\mathcal{W}^{N+1}$ is $\widetilde{\boldsymbol{F}}$-invariant subspace. Moreover, $\overline{\Psi}_{\lambda}(\mathcal{W})$ is submodule of $\overline{\Psi}_{\lambda}(\mathcal{V}).$
\end{prop}

From $\psi _0=\mathrm{id}_{\mathcal{V}^{N+1}},\, mc_0=\overline{\Psi}_0,$ the next proposition is obvious.

\begin{prop}\label{409}\ \ If $(\ast \ast)$ is satisfied, then $\overline{\Psi}_0(\mathcal{V}) \simeq \mathcal{V}.$
\end{prop}

Proof of the Proposition \ref{410}, \ref{411}, \ref{412} are similar to Detteweiler and Reiter's paper \cite{DR1}.

\begin{prop}\label{410}\ \ If $(\ast ),(\ast \ast)$ are satisfied, then for any $\lambda ,\mu \in \mathbb{C},\ \overline{\Psi}_{\mu} \circ \overline{\Psi}_{\lambda}(\mathcal{V}) \simeq \overline{\Psi}_{\mu}(\mathcal{V}^{N+1})/\overline{\Psi}_{\mu}(\mathcal{K}_{\mathcal{V}} +\mathcal{L}_{\mathcal{V}}(\lambda )).$
\end{prop}

$Proof.$\ \ If $\mu =0,$ it is easily seen that
\[ \overline{\Psi}_0 \circ \overline{\Psi}_{\lambda}(\mathcal{V}) \simeq \overline{\Psi}_{\lambda}(\mathcal{V}) = \mathcal{V}^{N+1}/(\mathcal{K}_{\mathcal{V}} +\mathcal{L}_{\mathcal{V}}(\lambda )) \simeq \overline{\Psi}_0(\mathcal{V}^{N+1} )/\overline{\Psi}_0(\mathcal{K}_{\mathcal{V}} +\mathcal{L}_{\mathcal{V}}(\lambda )).\]
Here we assume $\mu \ne 0.$ We set 
\begin{align}
\lambda '&=q^{\lambda}-1,\ \mu '=q^{\mu}-1,\ \mathcal{K}_1=\mathcal{K}_{\mathcal{V}},\ \mathcal{L}_1=\mathcal{L}_{\mathcal{V}}(\lambda ),\ \mathcal{K}_2=\mathcal{K}_{\mathcal{V}^{N+1}},\ \mathcal{L}_2=\mathcal{L}_{\mathcal{V}^{N+1}}(\mu ),\\
\widetilde{\boldsymbol{F}}&=\Psi _{\lambda }(\boldsymbol{B}),\ \check{\boldsymbol{F}}=\overline{\Psi}_{\lambda}(\boldsymbol{B}),\ \mathcal{M}=\overline{\Psi}_{\lambda}(\mathcal{V}),\ \mathcal{H}= \mathcal{K}_1+\mathcal{L}_1.
\end{align}
Let us first prove
\begin{equation}
\mathrm{(i)}\ \mathcal{K}_{\mathcal{M}}=(\mathcal{K}_2+\mathcal{H}^{N+1})/{\mathcal{H}^{N+1}},\ \ \ \mathrm{(ii)}\ \mathcal{L}_{\mathcal{M}}=(\mathcal{L}_2+\mathcal{H}^{N+1})/{\mathcal{H}^{N+1}}.
\end{equation}

(i)\ We set $\check{F}_0=1_m-\sum _{i=1}^N\check{F}_i-\check{F}_{\infty}.$ For any $k+\mathcal{H}^{N+1}={}^t({}^tk_0, \ldots ,{}^tk_N)+\mathcal{H}^{N+1} \in (\mathcal{K}_2+\mathcal{H}^{N+1})/{\mathcal{H}^{N+1}}$, we obtain $k+\mathcal{H}^{N+1} \in \mathcal{K}_{\mathcal{M}}$ from $\check{F}_i(k_i+\mathcal{H})=\mathcal{H}(i \in I=\{ 0, \ldots ,N \}).$ Therefore, $(\mathcal{K}_2+\mathcal{H}^{N+1})/{\mathcal{H}^{N+1}}$ is subspace of $\mathcal{K}_{\mathcal{M}}$.
On the other hand, for any $v+\mathcal{H}^{N+1}={}^t({}^tv_0, \ldots ,{}^tv_N)+\mathcal{H}^{N+1} \in \mathcal{K}_{\mathcal{M}}, v_i={}^t({}^tv_{i0}, \ldots ,{}^tv_{iN})\, (v_{ij}\in \mathcal{V}),$ we compute $\widetilde{F}_0v_0\, $: 
\[ \widetilde{F}_0v_0=(1_m-\textstyle{\sum _{i=1}^N}\widetilde{F}_i-\widetilde{F}_{\infty})v_0=(\widehat{F}-\textstyle{\sum _{i=1}^N}\widetilde{F}_i)v_0={}^t(\textstyle{\sum _{j=0}^N}{}^t(B_jv_{0j}),-\lambda '{}^tv_{01}, \ldots ,-\lambda '{}^tv_{0N})\]
and we find 
\[ \widetilde{F}_jv_j={}^t(0, \ldots ,\stackrel{j+1}{\stackrel{\vee}{\textstyle{\sum _{i=0}^N}{}^t(B_iv_{ji})+\lambda '{}^tv_{jj}}}, \ldots ,0)\, (j \in J=\{ 1, \ldots ,N \}).\]

(i-1)\ If $\lambda =0,$ then it is clear that $\Tilde{F}_iv_i={}^t(0, \ldots ,\sum _{j=0}^N{}^t(B_jv_{ij}), \ldots ,0)\, (i \in I).$ Moreover, $\widetilde{F}_iv_i \in \mathcal{H}= \mathcal{K}+\mathcal{L}= \mathcal{L}= \{ {}^t({}^tw_0, \ldots ,{}^tw_N) ; \sum _{j=0}^N B_jw_j=0 \} $ and $B_i\sum _{j=0}^NB_jv_{ij}=0.$ Hence we get
\[ \widetilde{F}_iv_i \in {}^t(0, \ldots ,\stackrel{i+1}{\stackrel{\vee}{\ker B_i}}, \ldots ,0).\]
Therefore, we obtain $v_i \in \ker \widetilde{F}_i+\mathcal{K}_1.$

(i-2)\ If $\lambda \ne 0,$ then
\[ \widetilde{F}_iv_i=({}^tk_{i0}+{}^th_i, \ldots ,{}^tk_{iN}+{}^th_i)\qquad (k_{ij} \! \in \! \ker B_j,\, h_i \in \ker (A_{\infty}-b_{\infty}q^{\lambda}1_m),\, i \in I).\]
If $i \ne 0,$ we get $h_i=-k_{ij} \in \ker B_j\, (j \in I\setminus \{ i\} ).$ Hence we see $h_i \in \ker (B_i+\lambda '1_m)$ from $h_i \in \ker (A_{\infty}-b_{\infty}q^{\lambda}1_m)= \ker (\sum _{r=0}^NB_r+\lambda '1_m).$ Since $(\ast \ast )$ is satisfied, we get $h_i=0$. Here
\[ \widetilde{F}_iv_i \in {}^t(0, \ldots ,\stackrel{i+1}{\stackrel{\vee}{\ker B_i}}, \ldots ,0).\]
If $i=0$, then it results in the case $i\ne 0$ because 
\begin{equation}
\widetilde{F}_0=1_{(N+1)m}-\sum _{r=1}^NF_r+\lambda '1_{(N+1)m}-F_{\infty}=
\begin{pmatrix} B_0+\lambda '1_m & \dotsi & B_N \\ & O & \end{pmatrix}.
\end{equation}

Hence we find $v+\mathcal{H}^{N+1} \in (\mathcal{K}_2+\mathcal{H}^{N+1})/{\mathcal{H}^{N+1}}.$ Moreover, $\mathcal{K}_{\mathcal{M}}$ is a subspace of $(\mathcal{K}_2+\mathcal{H}^{N+1})/{\mathcal{H}^{N+1}}$. Therefore, we obtain $\mathcal{K}_{\mathcal{M}}=(\mathcal{K}_2+\mathcal{H}^{N+1})/{\mathcal{H}^{N+1}}.$

(ii)\ For any
\[ v+\mathcal{H}^{N+1}\! ={}^t({}^th, \ldots ,{}^th)+\mathcal{H}^{N+1}\! \in \! (\mathcal{L}_2+\mathcal{H}^{N+1})/{\mathcal{H}^{N+1}}(h\! \in \! \ker (\widetilde{F}_{\infty}-q^{\mu}1_{(N+1)m})),\]
we let $\widetilde{H}=(F_{t-1})_{1 \leq s,t \leq N+1},\check{H}=(\check{F}_{t-1})_{1 \leq s,t \leq N+1}.$ Then we obtain
\[(\check{H}+\mu 1_{(N+1)^2m})(v+\mathcal{H}^{N+1})=(\widetilde{H}+\mu 1_{(N+1)^2m})v+\mathcal{H}^{N+1}=\mathcal{H}^{N+1}.\]
Consequently, we find $v+\mathcal{H}^{N+1} \in \mathcal{L}_{\mathcal{M}}$. Meanwhile, for any
\[ v+\mathcal{H}^{N+1}={}^t({}^th, \ldots ,{}^th)+\mathcal{H}^{N+1} \in \mathcal{L}_{\mathcal{M}}(h \in \ker (\overline{F}_{\infty}-q^{\mu}1_{(N+1)m})),\]
we see $v+\mathcal{H}^{N+1} \in (\mathcal{L}_2+\mathcal{H}^{N+1})/{\mathcal{H}^{N+1}}.$ Therefore, we obtain $\mathcal{L}_{\mathcal{M}}=(\mathcal{L}_2+\mathcal{H}^{N+1})/{\mathcal{H}^{N+1}}.$

Let us remember the isomorphism theorems. For a linear space $V$ and subspaces $W,W'$ of $V$,
\begin{align*}
&(\mathrm{iii})\ \mathrm{if}\ W' \subset W,\ \mathrm{then}\ (V/W')/(W/W') \simeq V/W;\\
&(\mathrm{iv})\ W'/(W \cap W') \simeq (W+W')/W.
\end{align*}
From the above, we can compute $mc_{\mu} \circ mc_{\lambda}(\mathcal{V})\, $:
\begin{align*}
mc_{\mu} \circ mc_{\lambda}(\mathcal{V}) &=mc_{\mu}(\mathcal{V}^{N+1}/{\mathcal{H}})\\
&=(\mathcal{V}^{N+1}/{\mathcal{H}})^{N+1}/(\mathcal{K}_{\mathcal{M}}+\mathcal{L}_{\mathcal{M}})\\
&=(\mathcal{V}^{{(N+1)}^2}/{\mathcal{H}}^{N+1})/((\mathcal{K}_2+\mathcal{L}_2+\mathcal{H}^{N+1})/\mathcal{H}^{N+1})\ \ (\because \ (\mathrm{i}),(\mathrm{ii}))\\
&\simeq (\mathcal{V}^{{(N+1)}^2}/(\mathcal{K}_2+\mathcal{L}_2))/((\mathcal{K}_2+\mathcal{L}_2+\mathcal{H}^{N+1})/(\mathcal{K}_2+\mathcal{L}_2))\ \ (\because \ (\mathrm{iii}))\\
&\simeq mc_{\mu}(\mathcal{V}^{N+1})/(\mathcal{H}^{N+1}/((\mathcal{K}_2+\mathcal{L}_2)\cap \mathcal{H}^{N+1}))\ \ (\because \ (\mathrm{iv}))\\
&=mc_{\mu}(\mathcal{V}^{N+1})/mc_{\mu}(\mathcal{K}_{\mathcal{V}} +\mathcal{L}_{\mathcal{V}}(\lambda )).\ \ \square 
\end{align*}

\begin{prop}\label{411}\ \ $mc_{\lambda}$ preserves conditions $(\ast ),(\ast \ast).$
\end{prop}

$Proof.$\ \ It is sufficient to prove that $\overline{\Psi}_{\lambda}$ preserves conditions $(\ast ),(\ast \ast).$ In the case $\lambda =0$ is obvious because of Proposition \ref{409}. Hence we assume $\lambda \ne 0$ and $\mathcal{V}$ satisfy $(\ast ),(\ast \ast).$ Here we use notations in proof of previous proposition. If $\tau =0,$ for any $v+\mathcal{H}={}^t({}^tv_0, \ldots ,{}^tv_N)+\mathcal{H}\in \bigcap _{i=0}^N \ker \check{F}_i,$ it is clear that $\widetilde{F}_0v \in \mathcal{H}.$ Here  we get $v \in \mathcal{H}$ from Proposition \ref{410}(i-2). Consequently, we obtain $\bigcap _{i=0}^N \ker \check{F}_i=\{ \mathcal{H}\} .$

If $\tau \ne 0,$ for any $v+\mathcal{H}\in \bigcap _{j \ne i} \ker \check{F}_j \cap (\check{F}_i+\tau 1_{(N+1)m})\, (i \in J=\{ 1, \ldots ,N \} )$, we get $v \in \mathcal{H}$ from $\widetilde{F}_0v \in \mathcal{H}.$ Hence we obtain $\bigcap _{j \ne i} \ker \check{F}_j \cap (\check{F}_i+\tau 1_{(N+1)m})=\{ \mathcal{H}\} .$ The case $i=0$ is reduced to the case $i\in J$. Therefore, $\overline{\Psi}_{\lambda}(\mathcal{V})$ satisfies $(\ast )$.

In the meantime, we put any $\tau \in \mathbb{C}$ and $v={}^t({}^tv_0, \ldots ,{}^tv_N) \in \mathcal{V}^{N+1}.$ If $i \in J,$ then 

\[ \widetilde{F}_iv={}^t(0, \ldots ,\stackrel{i+1}{\stackrel{\vee}{\sum _{j=0}^N{}^t(B_jv_j)+\lambda '{}^tv_i}} , \ldots ,0).\]

Hence $\widetilde{F}_iv$ spans the linear space ${}^t(0, \ldots ,\mathcal{V}, \ldots ,0).$ Moreover, it is clear that
\[ (\widetilde{F}_0+\tau 1_{(N+1)m})v={}^t(\sum _{j=0}^N{}^t(B_jv_j)+(\lambda '+\tau ){}^tv_0,\tau {}^tv_1, \ldots ,\tau {}^tv_N).\]
Consequently, $\sum _{j=0}^NB_jv_j+(\lambda '+\tau )v_0$ spans $\mathcal{V}$. Here the case $i=0$ is reduced to the case $i\in J$. Therefore, we obtain $\sum _{j \neq i}\mathrm{im}\check{F}_j+\mathrm{im}(\check{F}_i+\tau 1_{(N+1)m})=\mathcal{V}^{N+1}+\mathcal{H}\, (i \in J).$ From the above, $\overline{\Psi}_{\lambda}(\mathcal{V})$ satisfies $(\ast \ast )$.\ $\square$\\

Here the $\overline{\Psi}_{\lambda}$ satisfies the next proposition.

\begin{prop}\label{412}\ \ If $(\ast ),(\ast \ast)$ are satisfied, then for any $\lambda ,\mu \in \mathbb{C},\,  \overline{\Psi}_{\mu} \circ \overline{\Psi}_{\lambda}(\mathcal{V}) \simeq \overline{\Psi}_{\log _q(q^{\lambda}+q^{\mu}-1)}(\mathcal{V})$.
\end{prop}

$Proof.$\ \ If $\lambda \mu =0,$ it is obvious. We assume $\lambda \mu \ne 0$ and set
\begin{align}
&\widetilde{\boldsymbol{F}}=\Psi_{\lambda}(\boldsymbol{B}),\ \boldsymbol{F}'=\Psi_{\log _q(q^{\lambda}+q^{\mu}-1)}(\boldsymbol{B}),\ \boldsymbol{H}=\Psi_{\mu}(\widetilde{\boldsymbol{F}}),\ \mathcal{K}_1=\mathcal{K}_{\mathcal{V}},\ \mathcal{L}_1=\mathcal{L}_{\mathcal{V}}(\lambda ),\\
&\mathcal{K}_2=(\mathcal{K}_{\mathcal{V}^{N+1}},\widetilde{\boldsymbol{F}}),\ \mathcal{L}_2=(\mathcal{L}_{\mathcal{V}^{N+1}}(\mu ),\widetilde{\boldsymbol{F}}),\ \mathcal{L}'=\mathcal{L}_{\mathcal{V}}(\log _q(q^{\lambda}+q^{\mu}-1)),\ \mathcal{H}= \mathcal{K}_1+\mathcal{L}_1.
\end{align}
Here we prove that induced mapping $\overline{\phi}:\overline{\Psi}_{\mu} \circ \overline{\Psi}_{\lambda}(\mathcal{V}) \longrightarrow \overline{\Psi}_{\log _q(q^{\lambda}+q^{\mu}-1)}(\mathcal{V})$ is isomorphism from
\begin{equation}
\phi :\ \Psi_{\mu} \circ \Psi_{\lambda}(\mathcal{V}) \longrightarrow \Psi_{\log _q(q^{\lambda}+q^{\mu}-1)}(\mathcal{V})\ \biggr( {}^t({}^tv_0, \ldots ,{}^tv_N) \longmapsto \sum _{i=0}^N\widetilde{F}_iv_i \biggl).
\end{equation}
We first find 
\begin{equation}
\overline{\Psi}_{\mu} \circ \overline{\Psi}_{\lambda}(\mathcal{V}) \simeq \overline{\Psi}_{\mu}(\mathcal{V}^{N+1})/\overline{\Psi}_{\mu}(\mathcal{K}_{\mathcal{V}} +\mathcal{L}_{\mathcal{V}}(\lambda )) \simeq \mathcal{V}^{{(N+1)}^2}/(\mathcal{K}_2+\mathcal{L}_2+{\mathcal{H}}^{N+1}).
\end{equation}
It is easy to check that $(\mathcal{L}_1)^{N+1}$ is subspace of $\mathcal{K}_2=\ker (\phi )$. Moreover, we get $\phi ((\mathcal{K}_1)^{N+1})=\sum _{i=0}^N\widetilde{F}_i\mathcal{K}_1=\mathcal{K}_1$ and $\mathcal{L}_2=\{ {}^t({}^th, \ldots ,{}^th) ;h \in \ker (\widetilde{F}_{\infty}-q^{\mu}1_{(N+1)m})\}$. Hence we obtain 
\[ \phi (\mathcal{L}_2)=\sum _{i=0}^N\widetilde{F}_i\ker F_{\infty}'=\bigr( \sum _{i=0}^N\widetilde{F}_i\bigl) \ker F_{\infty}'=\ker F_{\infty}'=\mathcal{L}'\ (F_{\infty}' =\widetilde{F}_{\infty}-q^{\mu}1_{(N+1)m}).\]
Here we compute $\dim (\mathcal{K}_2)\, $:
\[ \dim (\mathcal{K}_2)=\sum _{i=0}^N\mathrm{dim \, ker}\widetilde{F}_i=\sum _{i=0}^N\{ \dim (\mathcal{V}^{N+1})-\mathrm{rank}\widetilde{F}_i \}=\sum _{i=0}^N\{ (N+1)m-m \}=N(N+1)m.\]
Consequently, we can calculate $\dim (\overline{\Psi}_{\mu} \circ \overline{\Psi}_{\lambda}(\mathcal{V}))\, $:
\begin{align*}
\dim (\overline{\Psi}_{\mu} \circ \overline{\Psi}_{\lambda}(\mathcal{V}))&=\dim (\mathcal{V}^{{(N+1)}^2}/(\mathcal{K}_2+\mathcal{L}_2+\mathcal{H}^{N+1}))\\
&=\dim (\mathcal{V}^{{(N+1)}^2})-\dim (\mathcal{K}_2+\mathcal{L}_2+(\mathcal{K}_1)^{N+1}+(\mathcal{L}_1)^{N+1})\\
&=(N+1)^2m-\dim (\mathcal{K}_2+\mathcal{L}_2+(\mathcal{K}_1)^{N+1})\\
&=(N+1)^2m-\dim (\mathcal{K}_2)-\dim (\mathcal{L}_2+(\mathcal{K}_1)^{N+1})\\
&=(N+1)^2m-N(N+1)m-\dim (\mathcal{K}_1+\mathcal{L}')\\
&=(N+1)m-\dim (\mathcal{K}_1+\mathcal{L}')\\
&=\dim (\mathcal{V}^{N+1})-\dim (\mathcal{K}_1+\mathcal{L}')\\
&=\dim (\mathcal{V}^{N+1}/(\mathcal{K}_1+\mathcal{L}')).
\end{align*}
Here we set $\lambda '=q^{\lambda}-1,\mu '=q^{\mu}-1$. For any
\[ v={}^t({}^tv_0, \ldots ,{}^tv_N) \in \mathcal{V}^{{(N+1)}^2},\ (v_j={}^t({}^tv_{j0}, \ldots ,{}^tv_{jN}),\, v_{ij}\in \mathcal{V}),\]
we get the following relations.
\begin{gather*}
(F_i'\circ \phi )(v)={}^t(0, \ldots ,\stackrel{i+1}{\stackrel{\vee}{{}^tw_i}} , \ldots , 0)=(\phi \circ H_i)(v)\, (i \in \{ 0, \ldots ,N \} ),\\
w_i=\sum _{j=0}^NB_j\biggr\{ \sum _{k=0}^NB_kv_{jk}+\lambda 'B_jv_{jj}+(\lambda '+\mu ')v_{ij} \biggl\}+\lambda '(\lambda '+\mu ')v_{ii},\\
F_{\infty}'\circ \phi =\phi -\sum _{i=0}^N(F_i'\circ \phi )=\phi -\sum _{i=0}^N(\phi \circ H_i )=\phi \circ H_{\infty}.
\end{gather*}
Therefore, we obtain $\overline{\Psi}_{\mu} \circ \overline{\Psi}_{\lambda}(\mathcal{V}) \simeq \overline{\Psi}_{\log _q(q^{\lambda}+q^{\mu}-1)}(\mathcal{V}).\ \square$\\

From the above, Theorem \ref{irre} is shown.\\
\\
\rm{\textbf{Theorem \ref{irre}}\ (irreducibility)}\ \ \it{If $(\ast ),(\ast \ast )$ are satisfied, then $\mathcal{V}$ is irreducible if and only if $mc_{\lambda}(\mathcal{V})$ is irreducible}\rm{.}\\

$Proof.$\ \ For any non-zero irreducible module $\mathcal{V}$ and $\lambda \in \mathbb{C}$, we put $\mathcal{M}=\overline{\Psi}_{\lambda}(\mathcal{V})$ and non-zero submodule $\mathcal{M}'$ of $\mathcal{M}$. Here $\mathcal{W}=\overline{\Psi}_{\log _q(1-q^{\lambda})}(\mathcal{M}')$ is submodule of
\begin{equation*}
\overline{\Psi}_{\log _q(1-q^{\lambda})}(\mathcal{M})=(\overline{\Psi}_{\log _q(1-q^{\lambda})}\circ \overline{\Psi}_{\lambda})(\mathcal{V}) \simeq \overline{\Psi}_0(\mathcal{V})=mc_0(\mathcal{V})\simeq \mathcal{V}.
\end{equation*}
Hence we obtain $\mathcal{W}=0\ \mathrm{or}\ \mathcal{V}$. If $\mathcal{W}=0,$ then we get $\mathcal{M}' \simeq \overline{\Psi}_{\lambda}(\mathcal{W})=\overline{\Psi}_{\lambda}(0)=0.$ This is a contradiction. Consequently, we find $\mathcal{W}=\mathcal{V}$. Moreover, we get
\begin{equation*}
\mathcal{M}'= \overline{\Psi}_{\lambda}(\mathcal{W})=\overline{\Psi}_{\lambda}(\mathcal{V})=\mathcal{M}.
\end{equation*}
Hence $\mathcal{W}$ is irreducible module. Here $\overline{\Psi}_{\lambda}(\mathcal{V})$ is irreducible if and only if $mc_{\lambda}(\mathcal{V})$ is irreducible. Therefore, $\mathcal{V}$ is irreducible if and only if $mc_{\lambda}(\mathcal{V})$ is irreducible. The proof of the theorem has been completed.\ $\square$\\

\subsection{Proof of Theorem \ref{index}.}

In this section, we prove that $mc_{\lambda}$ preserves rigidity index of equation $E_R$. At first, we examine the change of spectral types $S_0,S_{\infty}$.

\begin{lem}\ \ We set coefficient polynomial $A(x)=\sum _{k=0}^NA_kx^k\ (resp.\ G(x)=\sum _{k=0}^NG_kx^k)$ of canonical form of $E_{\boldsymbol{B},\boldsymbol{b}}\ (resp.\ E_{\boldsymbol{F},\boldsymbol{b}})$, we let  $\mathrm{Ev}(M)$ be the set of eigenvalues of $M\in \mathrm{M}_m(\mathbb{C})$. If
\[A_0\sim \bigoplus _{\theta \in \mathrm{Ev}(A_0)}\bigoplus _{j=1}^{s_{\theta}^0}J(\theta ,t_{\theta ,\, j}^0),\ \ \ A_{\infty}\sim \bigoplus _{\kappa \in \mathrm{Ev}(A_{\infty})}\bigoplus _{j=1}^{s_{\kappa}^{\infty}}J(\kappa ,t_{\kappa ,\, j}^{\infty})\]
and $(\ast \ast )$ is satisfied, then we obtain
\begin{align*}
G_0&\sim \bigoplus _{\theta \in \mathrm{Ev}(A_0)\setminus \{ q^{\lambda}\} }\bigoplus _{j=1}^{s_{\theta}^0}J(\theta ,t_{\theta ,\, j}^0)\oplus \bigoplus _{j=1}^{s_{q^{\lambda}}^0}J(q^{\lambda} ,t_{q^{\lambda} ,\, j}^0+1)\oplus J(q^{\lambda},1)^{\oplus (Nm-s_{q^{\lambda}}^0)},\\
G_{\infty}&\sim \bigoplus _{\kappa \in \mathrm{Ev}(A_{\infty})\setminus \{ b_{\infty}\} }\bigoplus _{j=1}^{s_{\kappa}^{\infty}}J(\kappa ,t_{\kappa ,\, j}^{\infty})\oplus \bigoplus _{j=1}^{s_{b_{\infty}}^{\infty}}J(b_{\infty} ,t_{b_{\infty} ,\, j}^{\infty}+1)\oplus J(b_{\infty},1)^{\oplus (Nm-s_{b_{\infty}}^{\infty})}.
\end{align*}
\end{lem}

$Proof.$\ \ It is easily seen that $G_0=1_m-F_0,\ G_{\infty}=b_{\infty}F_{\infty},\ F_0=1_m-\sum _{i=1}^NF_i-F_{\infty}$ and 
\begin{equation}
\theta 1_{(N+1)m}-G_0=\begin{pmatrix} \theta 1_m-A_0 & B_1& \dotsi & B_N\\
 & (\theta -q^{\lambda})1_m & & \\
 & & \ddots & \\
 & & & (\theta -q^{\lambda})1_m  \end{pmatrix}\ \ (\theta \in \mathbb{C}).
\end{equation}

(i)\ \ If $\theta \ne q^{\lambda},$ then $\dim \ker ((\theta 1_{(N+1)m}-G_0)^n)=\dim \ker ((\theta 1_m-A_0)^n)\ (n\in \mathbb{Z}_{>0}).$

(ii)\ \ If $\theta =q^{\lambda},$ then for any $v={}^t({}^tv_0,\, \ldots \, ,{}^tv_N)\in \mathcal{V}^{N+1}\ (v_i\in \mathcal{V}),$ we get
\[ (\theta 1_{(N+1)m}-G_0)v={}^t({}^tv',0,\, \ldots \, ,0),\, v'=\sum _{k=0}^NB_kv_k+(\theta -1)1_m.\]
Here $v'$ spans $\mathcal{V}$ because condition $(\ast \ast )$. Hence we obtain 
\[ \dim \ker (\theta 1_{(N+1)m}-G_0)=Nm,\ \dim \ker ((\theta 1_{(N+1)m}-G_0)^{n+1})=\dim \ker ((\theta 1_m-A_0)^n)\ (n\in \mathbb{Z}_{>0}).\]

(iii)\ \ If $\kappa \ne b_{\infty}$, then for any $v={}^t({}^tv_0,\, \ldots \, ,{}^tv_N)\in \ker ((\kappa 1_{(N+1)m}-G_{\infty})^n)\ (v_i\in \mathcal{V},\, n\in \mathbb{Z}_{>0})$, we get
\[ 0=(\kappa 1_{(N+1)m}-G_{\infty})^nv=\{ (\kappa -b_{\infty})1_{(N+1)m}+b_{\infty}\widehat F\} ^nv=(\kappa -b_{\infty})^nv+P\widehat{F}v\ (P\in \mathrm{M}_m(\mathbb{C})).\]
Hence we find $v_0=\cdots =v_N.$ Moreover, it is clear that
\[ (\kappa 1_{(N+1)m}-G_{\infty})^nv=\{ (\kappa -b_{\infty})1_{(N+1)m}+b_{\infty}\widehat{F}\} ^nv={}^t({}^tv',\, \ldots \, ,{}^tv'),\ v'=(\kappa 1_m-A_{\infty})^nv_0.\]
Therefore, we obtain $\dim \ker ((\kappa 1_{(N+1)m}-G_{\infty})^n)=\dim \ker ((\kappa 1_m-A_{\infty})^n)$.

(iv)\ \ If $\kappa =b_{\infty}$, then we obtain 
\[ \dim \ker (\kappa 1_{(N+1)m}-G_{\infty})=\dim \ker \widehat{F}=(N+1)m-\mathrm{dim\ im}\widehat{F}=(N+1)m-m=Nm\]
from $\kappa 1_{(N+1)m}-G_{\infty}=b_{\infty}\widehat F$ and $(\ast \ast ).$ Here for any
\[ v={}^t({}^tv_0,\, \ldots \, ,{}^tv_N)\in \ker ((\kappa 1_{(N+1)m}-G_{\infty})^{n+1})\ (v_i\in \mathcal{V},\, n\in \mathbb{Z}_{>0}),\]
it is easily seen that
\[ (\kappa 1_{(N+1)m}-G_{\infty})v=b_{\infty}\widehat{F}v={}^t({}^tv',\, \ldots \, ,{}^tv'),\ \ \ v'=b_{\infty}\sum _{k=0}^NB_kv_k\]
and 
\[ (\kappa 1_{(N+1)m}-G_{\infty})^{N+1}v={}^t({}^tw,\, \ldots \, ,{}^tw),\ \ \ w=(\kappa 1_m-A_{\infty})^nv'.\]
Therefore, we obtain $\dim \ker ((\kappa 1_{(N+1)m}-G_{\infty})^{n+1})=\dim \ker ((\kappa 1_m-A_{\infty})^n)$.\ $\square$\\

We prepare for examining changes of spectral type $S_{\mathrm{div}}$.

\begin{lem}\ \ We can reduce $G(x)$ to $\widetilde{G}(x):$
\begin{equation}
\widetilde{G}(x)=\begin{pmatrix} T(x)1_m & & & \\
 & \ddots & & \\
 & & T(x)1_m & \\
 V_1(x) & \cdots & V_N(x) & A(q^{-\lambda}x)  \end{pmatrix}
\end{equation}
by elementary matrices. Here $V_i(x)\, (i=1,\ldots ,N)$ are polynomials and $T(x)=\prod _{k=1}^N(1-\frac{x}{b_k}).$
\end{lem}

$Proof.$\ \ For any $\lambda \in \mathbb{C},k\in J=\{ 1,\, \ldots ,\, N\} ,b_k\in \mathbb{C}\setminus \{ 0\}$, let $s_k=1-\frac{x}{b_k},\ s_k'=1-\frac{x}{q^{\lambda}b_k},\ T_k=\frac{T(x)}{s_k},\ b_{i,j}=1-\frac{b_i}{b_j}$. It is clear that 
\begin{align}
G(x)&=T(x)F(x)\\
&=\left( \prod _{k=1}^Ns_k\right) \cdot \left( F_{\infty}+\sum _{l=1}^N\frac{F_l}{s_l}\right) \\
&=T(x)1_m\oplus \bigoplus _{k=1}^Nq^{\lambda}s_k'T_k(x)1_m+\left( -T(x)1_m\oplus \bigoplus _{k=1}^N\frac{xT_k(x)}{b_k}1_m\right) \begin{pmatrix} 1_m \\
 \vdots \\
 1_m  \end{pmatrix}(B_0\cdots B_n).
\end{align}
Here we row reduce $G(x)$ by the elementary matrix
\begin{equation}
\begin{pmatrix} (1-s_1)1_m & s_11_m & &\\
 -1_m & 1_m & & \\
 \vdots & & \ddots & \\
 -1_m & & & 1_m  \end{pmatrix}.
\end{equation}
Next, we column reduce by the elementary matrix
\begin{equation}
\begin{pmatrix} 1_m & & & \\
 1_m & 1_m & & \\
 \vdots & & \ddots & \\
 1_m & & & 1_m  \end{pmatrix}.
\end{equation}
Then we obtain
\begin{equation}
q^{\lambda}\begin{pmatrix} T1_m & s_1'T1_m& & \\
 & s_1'T_11_m & & \\
 & & \ddots & \\
 & & & s_N'T_N1_m  \end{pmatrix}+\begin{pmatrix} O_m \\
 T_11_m \\
 \vdots \\
 T_N1_m  \end{pmatrix}(q^{\lambda}1_m-B_{\infty} \ B_1\cdots B_n).
\end{equation}
We set
\begin{align}
f_{i,j}&=(-1)^{i+j}b_ib_j^{-1}b_{i+1,j}^{-1}\prod _{k=1}^{j-1}(b_{j,k}^{-1}b_{i,k})\cdot \prod _{k=j+1}^{i-1}(b_kb_j^{-1}b_{k,j}^{-1}b_{i,k})\ \ (b_{N+1,j}=1),\\
g_i&=-\prod _{k=1}^ib_{i+1,k}^{-1}\cdot \prod _{k=1}^{i-1}b_{i,k} (\ne 0)
\end{align}
and $C_0=(C_{i,j}^0)_{1\le i,j\le N+1}\in \mathrm{M}_{(N+1)m}(\mathbb{C})\ (C_{i,j}^1\in \mathrm{M}_m(\mathbb{C}))$ as
\begin{equation}
C_{i,j}^0=\begin{cases}
 1_m & (i=j=1)\\
 f_{i-1,j-1}s_{i-1}1_m & (2\le i,j\le N,i\ge j)\\
 g_{i-1}s_i1_m & (2\le i=j-1 \le N)\\
 f_{N,j-1}1_m & (i=N+1,j\ne 1)\\
 O_m & (\mathrm{otherwise})
\end{cases}.
\end{equation}
Here $C_0$ is an elementary matrix. Let $j\in \{ 1,\, \ldots ,\, N\},\, I_j=\{ 1,\, \ldots ,\, j\}$. We prove
\begin{align*}
\mathrm{(i)}&\ \sum _{k=1}^lf_{l,k}=-g_l\ \ (l\in I_{N-1}),\\
\mathrm{(ii)}&\ \sum _{k=1}^Nf_{N,k}T_k(x)=t_0\ \ \left( t_0=\prod _{k=1}^{N-1}b_{N,k}\ne 0\right) .
\end{align*}

(i)\ It is clear that $\sum _{k=1}^lf_{l,k}(b_{l+1})\equiv 0\ (\mathrm{mod}\, g_l(b_{l+1}))$. Here we set $f(b_{l+1})=-(g_l)^{-1}\sum _{k=1}^lf_{l,k},$ then we find $\deg f(b_{l+1})\le l-1$ and $f(b_s)=1\ (s\in I_l).$ Therefore, for any $b_{l+1}\in \mathbb{C}$, we obtain $f(b_{l+1})=1.$

(ii)\ Let $g(x)=\sum _{k=1}^Nf_{N,k}T_k(x),$ then we find $\deg g(x)\le N-1$ and $g(b_s)=t_0\ (s\in I_N).$ Therefore, $g(x)=t_0.$ Hence we get
\begin{equation}
C_0\begin{pmatrix} O_m \\
 T_11_m \\
 \vdots \\
 T_N1_m  \end{pmatrix}=\begin{pmatrix} O_m \\
 \vdots \\
 O_m \\
 t_01_m \end{pmatrix}.
\end{equation}
Here let us reduce
\begin{equation}
q^{\lambda}C_0\begin{pmatrix} T1_m & s_1'T1_m& & \\
 & s_1'T_11_m & & \\
 & & \ddots & \\
 & & & s_N'T_N1_m  \end{pmatrix}.
\end{equation}
We set $U_{i,j}(p):=(p\delta _{si}\delta _{tj}1_m)_{1\le s,t\le N+1}\in \mathrm{M}_{(N+1)m}(\mathbb{C})\ (p\in \mathbb{C})$ and
\begin{align}
h_{i,j}&=g_j^{-1}\sum _{k=1}^jf_{i,k},\\
C_l&=1_{(N+1)m}+\sum _{k=l+2}^NU_{k,l+1}(h_{k-1,l})\ (1\le l\le N-2),\\
C_{N-1}&=1_m\oplus (-g_1^{-1})1_m\oplus \cdots \oplus (-g_{N-1}^{-1})1_m \oplus 1_m,\\
C&=C_{N-1}C_{N-2}\cdots C_1C_0.
\end{align}
Then we obtain 
\begin{equation}
\begin{split}
&C\begin{pmatrix} T1_m & s_1'T1_m& & \\
 & s_1'T_11_m & & \\
 & & \ddots & \\
 & & & s_N'T_N1_m  \end{pmatrix}\\
&\ =\begin{pmatrix} T1_m & s_1'T1_m & \\
 & s_1'T1_m & -s_2'T1_m & & & \\
 & & s_2'T1_m & -s_3'T1_m & & \\
 & & & \ddots & \ddots & \\
 & & & & s_{N-1}'T1_m & -s_N'T1_m\\
 & f_{N,1}s_1'T_11_m & f_{N,2}s_2'T_21_m & \cdots & f_{N,N-1}s_{N-1}'T_{N-1}1_m & f_{N,N}s_N'T_N1_m  \end{pmatrix}.
\end{split}
\end{equation}
For any $i\in I_{N-1}$, we set
\begin{equation}
\begin{split}
u_i&=\prod _{k=1}^is_k',\ \ u_i'=\prod _{k=1}^ib_{i+1,k},\ \ \widetilde{u}_i=s_{i+1}^{-1}(u_i-u_i'),\\
D_0&=1_{(N+1)m}-U_{1,2}(s_1'),\\
D_{1,i}&=(1_{(N+1)m}+U_{i+2,i+1}(\widetilde{u}_i))(1_{(N+1)m}+U_{i+2,i+2}(u_i'-1))(1_{(N+1)m}+U_{i+1,i+2}(s_{i+1}')),\\
D_{2,i}&=(1_{(N+1)m}+U_{i+2,i+1}(-\widetilde{u}_is_{i+1}'))(1_{(N+1)m}+U_{i+1,i+1}(u_i'^{-1}-1))(1_{(N+1)m}+U_{i+1,i+2}(s_{i+1}')),\\
D_1&=D_0D_{1,1}\cdots D_{1,N-1},\ \ \ D_2=D_{2,N-1}\cdots D_{2,1}.
\end{split}
\end{equation}
Here we remember $A(q^{-\lambda}x)=T(q^{-\lambda}x)B(q^{-\lambda}x)=(\prod _{k=1}^Ns_k')(B_{\infty}+\sum _{l=1}^Ns_l'^{-1}B_l),$ and we compute 
\begin{equation}
D_2C\left\{ q^{\lambda}\begin{pmatrix} T1_m & s_1'T1_m& & \\
 & s_1'T_11_m & & \\
 & & \ddots & \\
 & & & s_N'T_N1_m  \end{pmatrix}+\begin{pmatrix} O_m \\
 T_11_m \\
 \vdots \\
 T_N1_m  \end{pmatrix}(q^{\lambda}1_m-B_{\infty} \ B_1\cdots B_n)\right\} D_1.
\end{equation}
This is the $\widetilde{G}(x)$.\ $\square$\\

We prove the next lemma for examining type of elementary divisors of $G(x)$.

\begin{lem}\ \ For coefficient polynomial $A(x)=\sum _{k=0}^NA_kx^k$ of canonical form of Fuchsian equation $E_R$, we define $P_A\in \mathrm{M}_{Nm}(C)$ as
\begin{equation}
\begin{pmatrix}  & 1_m & & \\
 & & \ddots & \\
 & & & 1_m \\
 -A_{\infty}^{-1}A_0 & -A_{\infty}^{-1}A_1 & \cdots & -A_{\infty}^{-1}A_{N-1}  \end{pmatrix}.
\end{equation}
Then, for any $a_i\in Z_R=\{a \in \mathbb{C}\, ;\, \mathrm{det}A(a)=0 \}$, we obtain
\begin{equation}
n_j^i=\dim \ker ((a_i 1_{Nm}-P_A)^j)-\dim \ker ((a_i 1_{Nm}-P_A)^{j-1})\, (j\in \mathbb{Z}_{>0}).
\end{equation}
\end{lem}

$Proof.$\ \ $x1_{Nm}-P_A$ can be transformed to $1_{(N-1)m}\oplus A(x)$ by elementary matrices. Therefore, type of elementary divisors of $x1_{Nm}-P_A$ and type of elementary divisors of $A(x)$ are equal except for $1_{(N-1)m}$.\ $\square$\\

We obtain the following lemma by calculating the dimensions of the generalized eigenspaces of $P_A$.

\begin{lem}\ \ Let $I_j=\{ 1,\, \ldots ,\, j\} ,j_1=\min \{ N+1,j\} ,j_2=\max \{ N+1,j\} \ (j\in \mathbb{Z}_{>0}),$
$I_j'=\{ 1,\, \ldots ,\, j_2\}$. For any $a\in \mathbb{C}\setminus \{ 0\}$, the following conditions are equivalent{\normalfont :}
\begin{align*}
(i)&\ {}^t({}^tv_1,\, \ldots \, ,{}^tv_N)\in \ker (a1_{Nm}-P_A)\ (v_k\in \mathcal{V}),\\
(ii)&\ There\ exist\ v_{j_1},\, \ldots ,\, v_{j_2}\in \mathcal{V}\ such\ that\ for\ w_k=\sum _{l=1}^k(-1)^{l-1}\left( \begin{matrix} k-1 \\ l-1 \end{matrix} \right) a^{k-l}v_l\ (k\in I_j),\\
&\ v_k=\sum _{l=1}^j(-1)^{l-1}\left( \begin{matrix} k-1 \\ l-1 \end{matrix} \right) a^{k-l}w_l\, (k\in I_j')\ and\ \sum _{i=0}^{k-1}\frac{(-1)^i}{i!}\frac{d^iA}{dx^i}(a)w_{i+j-k+1}=0\, (k\in I_j). 
\end{align*}
\end{lem}

$Proof.$\ \ If $j=1$, then for any $v={}^t({}^tv_1,\, \ldots \, ,{}^tv_N)\in \ker (a1_{Nm}-P_A)\ (v_k\in \mathcal{V}),$ we put $w_1=v_1,v_{N+1}=a^Nv_1.$ Here we get $v_k=a^{k-1}v_1=a^{k-1}w_1\, (k\in I_1')$ and $A(a)w_1=\sum _{k=0}^NA_ka^kv_1=\sum _{k=0}^NA_kv_{k+1}=0$ from $(P_A-a1_{Nm})v=0.$ We assume that the equivalence is satisfied in the case $j=j'\in \mathbb{Z}_{>0}$. For any $v={}^t({}^tv_1,\, \ldots \, ,{}^tv_N)\in \ker ((a1_{Nm}-P_A)^{j'+1})\ (v_k\in \mathcal{V}),$ we let
\[ u={}^t({}^tu_1,\, \ldots \, ,{}^tu_N)=(a1_{Nm}-P_A)v\, (u_k\in \mathcal{V}),\ v_{N+1}=-A_{\infty}^{-1}\sum _{k=0}^{N-1}A_kv_{k+1}.\]
Then we find $u_k=av_k-v_{k+1}$ and $\sum _{k=0}^NA_kv_{k+1}=0\ (k\in I_0=\{ 1,\, \ldots ,\, N\})$. Here we set 
\[ \widetilde{w}_k=\sum _{l=1}^k(-1)^{l-1}\! \left( \begin{matrix} k-1 \\ l-1 \end{matrix} \right) \! a^{k-l}u_l\ (k\in I_{j'}) .\]
There exist $u_{j_1'},\, \ldots ,\, u_{j_2'}\in \mathcal{V}\ (j_1'=\min \{ N+1,j'\} ,j_2'=\max \{ N+1,j'\} )$ such that 
\begin{equation}
\widetilde{w}_k=\sum _{l=1}^k(-1)^{l-1}\left( \begin{matrix} k-1 \\ l-1 \end{matrix} \right) a^{k-l}u_l=\sum _{l=1}^{k+1}(-1)^{l-1}\left( \begin{matrix} k \\ l-1 \end{matrix} \right) a^{k+1-l}v_l\ \ (k\in I_{j'}).
\end{equation}
Let $w_1=v_1,w_k=\widetilde{w}_{k-1}\, (k\in I_{j'+1}\setminus \{ 1\} ),$ we obtain
\[ w_k=\sum _{l=1}^k(-1)^{l-1}\left( \begin{matrix} k-1 \\ l-1 \end{matrix} \right) a^{k-l}v_l\ \ (k\in I_{j'+1}).\]
Here we find $\displaystyle{u_k=\sum _{l=1}^{j'}(-1)^{l-1}\left( \begin{matrix} k-1 \\ l-1 \end{matrix} \right) a^{k-l}w_{l+1}\, (k\in I_{j'}').}$ We put $v_k\in \mathcal{V}$ such that $av_k-v_{k+1}=u_k\ (k\in \{ j_1',\, \ldots ,\, j_2'\})$. For any $k\in I_{j'}'$, we get
\begin{equation}
av_k-v_{k+1}=\sum _{l=1}^{j'}(-1)^{l-1}\left( \begin{matrix} k-1 \\ l-1 \end{matrix} \right) a^{k-l}w_{l+1}=\sum _{l=2}^{j'+1}(-1)^{l-2}\left( \begin{matrix} k-1 \\ l-2 \end{matrix} \right) a^{k+1-l}w_l
\end{equation}
and
\begin{equation}
v_k=\sum _{l=1}^{j'+1}(-1)^{l-1}\left( \begin{matrix} k-1 \\ l-1 \end{matrix} \right) a^{k-l}w_l.
\end{equation}
Moreover, we obtain 
\[ \sum _{i=0}^{k-1}\frac{(-1)^i}{i!}\frac{d^iA}{dx^i}(a)w_{i+(j'+1)-k+1}=0\]
from $\displaystyle{\sum _{i=0}^{k-1}\frac{(-1)^i}{i!}\frac{d^iA}{dx^i}(a)\widetilde{w}_{i+j'-k+1}=0\ \ (k\in I_{j'}).}$
On the other hand, by the computation:
\begin{equation}
\begin{split}
0&=\sum _{k=0}^NA_kv_{k+1}\\
&=\sum _{k=0}^{N-1}A_k\sum _{l=1}^{j'+1}(-1)^{l-1}\left( \begin{matrix} k \\ l-1 \end{matrix} \right) a^{k+1-l}w_l\\
&\ \ \ \ +A_N\left\{ \sum _{l=1}^{j'+1}(-1)^{l-1}\left( \begin{matrix} N-1 \\ l-1 \end{matrix} \right) a^{N+1-l}w_l-\sum _{l=1}^{j'}(-1)^{l-1}\left( \begin{matrix} N-1 \\ l-1 \end{matrix} \right) a^{N-l}w_{l+1}\right\} \\
&=\sum _{l=0}^{(j'+1)-1}\frac{(-1)^l}{l!}\sum _{k=0}^N\frac{k!}{(k-l)!}A_k a^{k-l}w_{l+1}\\
&=\sum _{l=0}^{(j'+1)-1}\frac{(-1)^l}{l!}\frac{d^lA}{dx^l}(a)w_{l+(j'+1)-(j'+1)+1},
\end{split}
\end{equation}

(ii) is satisfied in the case $j=j'+1\in \mathbb{Z}_{>0}$. The proof of the lemma has been completed.\ $\square$\\

From the above, we can calculate the type of elementary divisors of $G(x)=c_{\lambda}(A)(x).$ We obtain the next lemma by calculating the dimension of the generalized eigenspaces of $P_{\widetilde{G}}\in \mathrm{M}_{N(N+1)m}$.

\begin{lem}\ \ If $(\ast) ,(\ast \ast )$ are satisfied, then for any $a\in Z_R=\{a \in \mathbb{C}\, ;\, \mathrm{det}A(a)=0 \}$ and $j\in \mathbb{Z}_{>0},$ we obtain {\rm (i),(ii):}\\ 
{\rm (i)}\ \ If $q^{\lambda}a\in q^{\lambda}Z_R\setminus \{ b_k\, ;\, k\in \{ 1,\, \ldots ,\, N\} \}$,\ then\\
 $\dim \ker ((q^{\lambda}a1_{N(N+1)m}-P_{\widetilde{G}})^j)=\dim \ker ((a1_{Nm}-P_A)^j).$\\
{\rm (ii)}\ \ If $q^{\lambda}a\in q^{\lambda}Z_R\cap \{ b_k\, ;\, k\in \{ 1,\, \ldots ,\, N\} \}$,\ then\ $\dim \ker (q^{\lambda}a1_{N(N+1)m}-P_{\widetilde{G}})=Nm$\ and\\
$\dim \ker ((q^{\lambda}a1_{N(N+1)m}-P_{\widetilde{G}})^{j+1})=\dim \ker ((a1_{Nm}-P_A)^j).$
\end{lem}

$Proof.$\ \ (i)\ For any
\[ v={}^t({}^tv_1,\, \ldots \, ,{}^tv_N)\in \ker (a1_{N(N+1)m}-P_{\widetilde{G}}),\ (v_k={}^t({}^tv_{k,0},\, \ldots \, ,{}^tv_{k,N}),\ v_{k,l}\in \mathcal{V}),\]
we find $v_k=a^{k-1}v_1,\widetilde{G}(q^{\lambda}a)v_1=0$. Moreover, we obtain $A(a)v_{1,N}=0$, $\dim \ker (q^{\lambda}a1_{N(N+1)m}-P_{\widetilde{G}})=\dim \ker A(a)=\dim \ker (a1_{Nm}-P_A)$ from $\widetilde{G}(q^{\lambda}a)v_1=0\ \Leftrightarrow \ v_{1,j}=0\, (j\ne N)$. Meanwhile, we assume $\dim \ker ((q^{\lambda}a1_{N(N+1)m}-P_{\widetilde{G}})^{j'})=\dim \ker ((a1_{Nm}-P_A)^{j'})\, (j=j'\in \mathbb{Z}_{>0})$. In another expression, for $w_k={}^t({}^tw_{k,0},\, \ldots \, ,{}^tw_{k,l})\in \mathcal{V}^{N+1}\, (w_{k,N}\in \mathcal{V},\, k\in J=\{ 1,\, \ldots ,\, j'\}),$ 
\begin{equation}
\begin{split}
&\sum _{i=0}^{k-1}\frac{(-1)^i}{i!}\frac{d^i\widetilde{G}}{dx^i}(q^{\lambda}a)w_{i+j'-k+1}=0\\
&\ \Leftrightarrow \ \sum _{i=0}^{k-1}q^{(j-i-1)\lambda}\frac{(-1)^i}{i!}\frac{d^iA}{dx^i}(a)w_{i+j'-k+1,N}=0,\ w_{k,l}=0\, (l\ne N). 
\end{split}
\end{equation}
Here if there exist
\[ w_k={}^t({}^tw_{k,0},\, \ldots \, ,{}^tw_{k,N})\in \mathcal{V}^{N+1}\, (w_{k,l}\in \mathcal{V},\, k\in J'=\{ 1,\, \ldots ,\, j'+1\})\]
such that $\sum _{i=0}^{k-1}\frac{(-1)^i}{i!}\frac{d^i\widetilde{G}}{dx^i}(q^{\lambda}a)w_{i+j'-k+2}=0.$ Then we get $w_{k,l}=0\, (k\ne 1,l\ne N).$ Moreover, we find
\[ w_{1,l}=0,\ \ \sum _{i=0}^{k-1}q^{(j'-i)\lambda}\frac{(-1)^i}{i!}\frac{d^iA}{dx^i}(a)w_{i+j'-k+2,N}=0\ (k\in J',\, l\ne N),\]
because $\sum _{i=0}^{j'}\frac{(-1)^i}{i!}\frac{d^i\widetilde{G}}{dx^i}(q^{\lambda}a)w_{i+1}=0.$ Therefore, we obtain 
\[ \dim \ker ((q^{\lambda}a1_{N(N+1)m}-P_{\widetilde{G}})^{j'+1})=\dim \ker ((a1_{Nm}-P_A)^{j'+1}).\]

(ii)\ If $q^{\lambda}a=b_{k_0}\in q^{\lambda}Z_R\cap \{ b_k\, ;\, k\in \{ 1,\, \ldots ,\, N\} \} \, (k_0\in \{ 1,\, \ldots ,\, N\}),$ then we obtain 
\[ \dim \ker (q^{\lambda}a1_{N(N+1)m}-P_{\widetilde{G}})=\dim \ker \, \widetilde{G}(k_0)=\dim \ker \, G(k_0)=(N+1)m-\mathrm{dim\ im}G(k_0)=Nm.\]
We assume that there exist $w_k={}^t({}^tw_{k,0},\, \ldots \, ,{}^tw_{k,N})\in \mathcal{V}^{N+1}\, (w_{k,l}\in \mathcal{V},\, k=1,2)$ such that 
\[ \widetilde{G}(q^{\lambda}a)w_2=0,\ \ \ \frac{d\widetilde{G}}{dx}(q^{\lambda}a)w_2=\widetilde{G}(q^{\lambda}a)w_1.\]
Then it is clear that $\frac{dT}{dx}(b_{k_0})\ne 0$. Hence we get 
\[ w_{2,l}=0\, (l\ne N),\ \ \ A(a)w_{2,N}=0,\ \ \ \frac{dA}{dx}(a)w_{2,N}=q^{\lambda}\sum _{l=0}^NU_l(q^{\lambda}a)w_{1,l}.\]
Here $q^{\lambda}\sum _{l=0}^NU_l(q^{\lambda}a)w_{1,l}$ spans $\mathcal{V}$ from condition $(\ast \ast )$.\ Moreover, we find
\[ \dim \ker ((q^{\lambda}a1_{N(N+1)m}-P_{\widetilde{G}})^2)=\dim \ker A(a)=\dim \ker (a1_{Nm}-P_A).\]
Therefore, we obtain
\[ \dim \ker ((q^{\lambda}a1_{N(N+1)m}-P_{\widetilde{G}})^{j'+2})=\dim \ker ((a1_{Nm}-P_A)^{j'+1}).\ \square\]

From the above, the next proposition is obvious.

\begin{prop}
If $(\ast) ,(\ast \ast )$\ are satisfied and the spectral type $S(E_R)=(S_0;S_{\infty};S_{\mathrm{div}})$ of Fuchsian equation $E_R$ is given as
\begin{equation}
\begin{split}
S_{\xi}&\, :\, m_{1,1}^{\xi} \ldots m_{1,{t_{1,1}^{\xi}}}^{\xi}, \ldots ,m_{{l_{\xi}},1}^{\xi} \ldots m_{{l_{\xi}},{t_{{l_{\xi}},1}^{\xi}}}^{\xi}\ (\xi =0,\infty ),\\
S_{\mathrm{div}}&\, :\, n_1^1 \ldots n_{k_1}^1, \ldots ,n_1^l \ldots n_{k_l}^l,
\end{split}
\end{equation}
then spectral type $S(c_{\lambda}(E_R))=(S_0';S_{\infty}';S_{\mathrm{div}}')$ satisfies
\begin{equation}
\begin{split}
S_0'&\, :\, \begin{cases}
 Nm\ m_{1,1}^0 \ldots m_{1,{t_{1,1}^0}}^0, \ldots ,m_{l_0,1}^0 \ldots m_{{l_0},{t_{l_0,1}^0}}^0 & (q^{\lambda}=\alpha _1^0)\\
 Nm,m_{1,1}^0 \ldots m_{1,{t_{1,1}^0}}^0, \ldots ,m_{l_0,1}^0 \ldots m_{{l_0},{t_{l_0,1}^0}}^0 & (q^{\lambda}\notin \mathrm{Ev}(A_0))
\end{cases},\\
S_{\infty}'&\, :\, \begin{cases}
 Nm\ m_{1,1}^{\infty} \ldots m_{1,{t_{1,1}^{\infty}}}^{\infty}, \ldots ,m_{{l_{\infty}},1}^{\infty} \ldots m_{{l_{\infty}},{t_{{l_{\infty}},1}^{\infty}}}^{\infty} & (b_{\infty}=\alpha _1^{\infty})\\
 Nm,m_{1,1}^{\infty} \ldots m_{1,{t_{1,1}^{\infty}}}^{\infty}, \ldots ,m_{{l_{\infty}},1}^{\infty} \ldots m_{{l_{\infty}},{t_{{l_{\infty}},1}^{\infty}}}^{\infty} & (b_{\infty}\notin \mathrm{Ev}(A_{\infty}))
\end{cases},\\
S_{\mathrm{div}}'&\, :\, \underbrace{Nm,\ldots ,Nm}_{r_1},Nm\ n_1^1 \ldots n_{k_1}^1, \ldots ,Nm\ n_1^{r_2} \ldots n_{k_{r_2}}^{r_2},n_1^{r_2+1} \ldots n_{k_{r_2+1}}^{r_2+1} ,\ldots ,n_1^l \ldots n_{k_l}^l\\
&\ \ \ \ (b_1,\ldots ,b_{r_1}\in \{ b_k\, ;\, k\in \{ 1,\ldots ,N\} \} \setminus q^{\lambda}Z_A,\ q^{\lambda}a_1,\ldots ,q^{\lambda}a_{r_2}\in \{ b_k\, ;\, k\in \{ 1,\ldots ,N\} \} ).
\end{split}
\end{equation}
\end{prop}

We show the next lemma in order to examine how $q$-middle convolution changes the spectral type.

\begin{lem}\ \ If $\lambda \ne 0, $ for $\theta ,\kappa ,a\in \mathbb{C}\setminus \{ 0\} ,I=\{ 1,\, \ldots ,\, N\} ,$ we obtain
\begin{align}
&\dim (\ker (\theta 1_{(N+1)m}-G_0) \cap \mathcal{K})=\begin{cases}
 \dim \ker (A_0-1_m) & (\theta =1)\\
 \sum _{k=1}^N\dim \ker B_k & (\theta =q^{\lambda})\\
 0 & (\theta \ne 1,q^{\lambda})
\end{cases},\\
&\dim (\ker (\kappa 1_{(N+1)m}-G_{\infty}) \cap \mathcal{K})=\begin{cases}
 \dim \ker (A_0-1_m)+\sum _{k=1}^N\dim \ker B_k & (\kappa =b_{\infty})\\
 0 & (\kappa \ne b_{\infty})
\end{cases},\\
&\dim (\ker \, G(a)\cap \mathcal{K})\\
&\ \ \ \ \ \ \ \ =\begin{cases}
 \dim \ker (A_0-1_m)+\sum _{k\ne j}\dim \ker B_k & (a=b_j)\\
 \dim \ker B_j & (a=q^{\lambda}b_j\in q^{\lambda}Z_A\setminus \{ b_k;k\in I\})\\
 0 & (\mathrm{otherwise})
\end{cases},\\
&\dim \left( \frac{dG}{dx}(a)\frac{}{}^{-1}(\mathrm{im}\, G(a))\cap \ker \, G(a)\cap \mathcal{K}\right) =\begin{cases}
 \dim \ker B_j & (a=q^{\lambda}b_j\in q^{\lambda}Z_A)\\
 0 & (\mathrm{otherwise})
\end{cases},\\
&\dim (\ker (\theta 1_{(N+1)m}-G_0) \cap \mathcal{L})=\begin{cases}
 \dim \ker (A_{\infty}-q^{\lambda}b_{\infty}1_m) & (\theta =q^{\lambda})\\
 0 & (\theta \ne q^{\lambda})
\end{cases},
\end{align}
\begin{align}
&\dim (\ker (\kappa 1_{(N+1)m}-G_{\infty}) \cap \mathcal{L})=\begin{cases}
 \dim \ker (A_{\infty}-q^{\lambda}b_{\infty}1_m) & (\kappa =q^{\lambda}b_{\infty})\\
 0 & (\kappa \ne q^{\lambda}b_{\infty})
\end{cases},\\
&\dim (\ker \, G(a)\cap \mathcal{L})=\begin{cases}
 \dim \ker (A_{\infty}-q^{\lambda}b_{\infty}1_m) & (a\in \{ b_k;k\in I\} )\\
 0 & (a\notin \{ b_k;k\in I\} )
\end{cases}.
\end{align}
\end{lem}

$Proof.$\ \ (i)\ (Change of $S_0$ due to the $\mathcal{K}$)
For $\theta \in \mathbb{C}$ and any $v={}^t({}^tv_0, \ldots ,{}^tv_N) \in \ker (\theta 1_{(N+1)m}-G_0) \cap \mathcal{K}\, (v_k\in \mathcal{V}),$ it is easily seen that 
\[ 0=(\theta 1_{(N+1)m}-G_0)v={}^t(\textstyle{\sum _{k=0}^N}{}^t(B_kv_k)+(\theta -1)\, {}^tv_0,(\theta -q^{\lambda})\, {}^tv_1 , \ldots ,(\theta -q^{\lambda})\, {}^tv_N).\]
If $\theta =1,$ then it is clear that $\theta \ne q^{\lambda}$ and $v_k=0\, (k\in I=\{ 1,\, \ldots ,\, N\} ),v_0 \in \ker (A_0-1_m).$ Here we get $\dim (\ker (\theta 1_{(N+1)m}-G_0) \cap \mathcal{K})=\dim \ker (A_0-1_m).$

If $\theta =q^{\lambda}, $ then we find $v_k\in \ker B_k\, (k\in I).$ Therefore, we obtain $\dim (\ker (\theta 1_{(N+1)m}-G_0) \cap \mathcal{K})=\sum _{k=1}^N\dim \ker B_k.$

(ii)\ (Change of $S_{\infty}$ due to the $\mathcal{K}$)
For $\kappa \in \mathbb{C}$ and any $v \in \ker (\kappa 1_{(N+1)m}-G_{\infty}) \cap \mathcal{K},$ we get
\[ 0=(\kappa 1_{(N+1)m}-G_{\infty})v=(\kappa 1_{(N+1)m}-b_{\infty}F_{\infty})v=\{ \kappa 1_{(N+1)m}-b_{\infty}(1_{(N+1)m}-\widehat{F})\} v=(\kappa -b_{\infty})v.\]
If $\kappa =b_{\infty},$ then we obtain $\dim (\ker (\kappa 1_{(N+1)m}-G_{\infty}) \cap \mathcal{K})=\dim \mathcal{K}=\dim \ker (A_0-1_m)+\sum _{k=1}^N\dim \ker B_k$.

(iii)\ (Change of $S_{\mathrm{div}}$ due to the $\mathcal{K}$)

(iii-a)
For any $v={}^t({}^tv_0, \ldots ,{}^tv_N) \in \ker \, G(b_k)\cap \mathcal{K}\, (v_k\in \mathcal{V},\, k \in I),$ it is clear that $v_k=0.$ Hence we get $\dim (\ker \, G(b_k)\cap \mathcal{K})=\dim \ker (A_0-1_m)+\sum _{l\ne k}\dim \ker B_l.$

(iii-b)
If $q^{\lambda}a_i\in q^{\lambda}Z_A\setminus \{ b_k\, ;\, k\in I\}$, then $T(q^{\lambda}a_i)\ne 0$. Hence we obtain 
\begin{equation}
\ker \, G(q^{\lambda}a_i)=\ker F(q^{\lambda}a_i)=\ker \left( 1_{(N+1)m}-\widehat{F}+\sum _{k=1}^N\frac{F_k}{1-\frac{q^{\lambda}a_i}{b_k}} \right) .
\end{equation}
For any $v={}^t({}^tv_0, \ldots ,{}^tv_N) \in \ker G(q^{\lambda}a_i)\cap \mathcal{K}\, (v_k\in \mathcal{V})$, we get 
\[ 0=\{ 1-\widehat{F}+\sum _{k=1}^N(1-q^{\lambda}a_ib_k^{-1})^{-1}F_k\} v=\{ 1_m\oplus _{k=1}^Nq^{\lambda}(1-a_ib_k^{-1})(1-q^{\lambda}a_ib_k^{-1})^{-1}1_m\} v.\] 
Here if $a_i\notin \{ b_k\, ;\, k\in I\}$, then $v=0.$ In the meantime, if $a_i=b_j\, (j\in I)$ and $k\ne j$, then $v_k=0$. Therefore, we find $v_j\in \ker B_j.$ From the above, we obatin 
\begin{equation}
\dim \left( \frac{dG}{dx}(q^{\lambda}a_i)\frac{}{}^{-1}(\mathrm{im}\, G(q^{\lambda}a_i))\cap \ker \, G(q^{\lambda}a_i)\cap \mathcal{K}\right) =\begin{cases}
 0 & (a_i\notin \{ b_k\, ;\, k\in I\} )\\
 \dim \ker B_j=n_1^j & (a_i=b_j)
\end{cases}.
\end{equation}

(iii-c)
If $q^{\lambda}a_i=b_{j'}\in q^{\lambda}Z_A\cap \{ b_k\, ;\, k\in I\} \, (j'\in I),$ then we put $w_k={}^t({}^tw_{k,0}, \ldots ,{}^tw_{k,N}) \in \mathcal{V}^{N+1}\, (w_k\in \mathcal{V},\, k=1,2)$ such that
\[ w_2\in \ker G(q^{\lambda}a_i)\cap \mathcal{K},\ \ \ \frac{dG}{dx}(q^{\lambda}a_i)w_2=G(q^{\lambda}a_i)w_1.\]
Hence we find $\ker G(q^{\lambda}a_i)\cap \mathcal{K}=\ker F_{j'}\cap \mathcal{K}$ and $q^{\lambda}\ne 1.$ Therefore, we get $w_{2,j'}=0$. Moreover, $G(q^{\lambda}a_i)w_1$ spans ${}^t(0,\ldots ,0,\mathcal{V},0,\ldots ,0)$ from $(\ast \ast )$.

If $a_i\notin \{ b_k\, ;\, k\in I\}$, then we get $w_{2,k}=0\, (k\ne j')$ from $\frac{dG}{dx}(q^{\lambda}a_i)w_2=G(q^{\lambda}a_i)w_1$. Therefore, $w_2=0.$ Meanwhile, if $a_i=b_j\, (j\in I)$ and $k\ne j,$ then we find $w_{2,k}=0$ and $w_{2,j}\in \ker B_j.$ From the above, we obatin 
\begin{equation}
\dim \left( \frac{dG}{dx}(q^{\lambda}a_i)\frac{}{}^{-1}(\mathrm{im}\, G(q^{\lambda}a_i))\cap \ker \, G(q^{\lambda}a_i)\cap \mathcal{K}\right) =\begin{cases}
 0 & (a_i\notin \{ b_k\, ;\, k\in I\} )\\
 \dim \ker B_j & (a_i=b_j)
\end{cases}.
\end{equation}

(iv)\ (Change of $S_0$ due to the $\mathcal{L}$)
For $\theta \in \mathbb{C}$ and any $v={}^t({}^th, \ldots ,{}^th) \in \ker (\theta 1_{(N+1)m}-G_0) \cap \mathcal{L}\ (h \in \ker (A_{\infty}-q^{\lambda}b_{\infty}1_m))$, we find 
\[ 0=(\theta 1_{(N+1)m}-G_0)v={}^t((\theta -q^{\lambda})\, {}^th , \ldots ,(\theta -q^{\lambda})\, {}^th).\]
If $\theta =q^{\lambda},$ then $h \in \ker (A_{\infty}-q^{\lambda}b_{\infty}1_m)$. Therefore, we obtain $\dim (\ker (\theta 1_{(N+1)m}-G_0) \cap \mathcal{L})=\dim \mathcal{L}=\dim \ker (A_{\infty}-q^{\lambda}b_{\infty}1_m).$

(v)\ (Change of $S_{\infty}$ due to the $\mathcal{L}$)
For $\kappa \in \mathbb{C}$ and $v={}^t({}^th, \ldots ,{}^th) \in \ker (\kappa 1_{(N+1)m}-G_{\infty}) \cap \mathcal{L}\ (h \in \ker (A_{\infty}-q^{\lambda}b_{\infty}1_m)),$ we get 
\[ 0=(\kappa 1_{(N+1)m}-G_{\infty})v=(\kappa -q^{\lambda}b_{\infty})v.\]
If $\kappa =q^{\lambda}b_{\infty},$ then we obtain $\dim (\ker (\kappa 1_{(N+1)m}-G_{\infty}) \cap \mathcal{L})=\dim \mathcal{L}=\dim \ker (A_{\infty}-q^{\lambda}b_{\infty}1_m)$.

(vi)\ (Change of $S_{\mathrm{div}}$ due to the $\mathcal{K}$)
For any $k \in I$, $\mathcal{L}$ is subspace of $\ker G(b_k)=\ker F_k.$ Therefore, we obtain 
\[ \dim (\ker G(b_k) \cap \mathcal{L})=\dim \mathcal{L}=\dim \ker (A_{\infty}-q^{\lambda}b_{\infty}1_m).\ \square \]

From the above, Theorem \ref{index} is shown.\\
\\
\rm{\textbf{Theorem \ref{index}}\ (rigidity index)}\ \ \it{If $(\ast ) ,(\ast \ast )$ are satisfied, then $mc_{\lambda}$ preserves rigidity index of Fuchsian equation $E_R$}\rm{.}\\

$Proof.$\ \ In the case $\lambda =0,$ it is obvious from Proposition \ref{404}. We assume $\lambda \ne 0$. Let coefficient $\overline{G}(x)=\sum _{k=0}^N\overline{G}_kx^k\ (\overline{G}_{\infty}=\overline{G}_N)$ of canonical form of $E_{\boldsymbol{\overline{F}},\boldsymbol{b}}\, (\boldsymbol{\overline{F}}=mc_{\lambda}(\boldsymbol{B}))$. It is clear that $q^{\lambda}\ne 1,q^{\lambda}b_{\infty}\ne b_{\infty}.$ Here let $\alpha _{i_0}^0=1,\alpha _{i_{\infty}}^{\infty}=q^{\lambda}b_{\infty}$. we get 
\begin{equation}
\dim \ker (A_0-1_m)=m_{i_0,1}^0,\ \ \dim \ker (A_{\infty}-q^{\lambda}b_{\infty}1_m)=m_{i_{\infty},1}^{\infty}.
\end{equation}
Moreover, we set 
\begin{equation}
b_k=\begin{cases}
 a_k & (k\in \{ 1,\ldots ,r\} )\\
 c_k & (k\in \{ r+1,\ldots ,N\} ,c_k\notin Z_A)\end{cases},\ \ d_k=\dim \ker B_k,\ \ d=\sum _{k=1}^Nd_k .
\end{equation}
Then we find
\begin{equation}
\dim (mc_{\lambda}(\mathcal{V}))=(N+1)m-m_{i_0,1}^0-m_{i_{\infty},1}^{\infty}-d,\ \ \ d=\sum _{k=1}^rn_1^k.
\end{equation}
Since these relations, we obtain 
\begin{align}
p_0&=\dim \ker (\overline{G}_0-q^{\lambda}1_{\dim (mc_{\lambda}(\mathcal{V}))})=Nm-m_{i_{\infty},1}^{\infty}-d,\\
p_{\infty}&=\dim \ker (\overline{G}_{\infty}-b_{\infty}1_{\dim (mc_{\lambda}(\mathcal{V}))})=Nm-m_{i_0,1}^0-d,\\
p_k&=\dim \ker \overline{G}(b_k)=Nm-m_{i_0,1}^0-m_{i_{\infty},1}^{\infty}-d+d_k\ (k\in \{ 1,\ldots ,N\} ).
\end{align}
From the above, rigidity index, $\mathop{\rm idx} (mc_{\lambda}(E_R))$, of equation $E_R$ is calculated:
\begin{align*}
\mathop{\rm idx} &(mc_{\lambda}(E_R))\\
&=\sum _{i \ne i_0}\sum _{j=1}^{t_{i,1}^0}(m_{i,j}^0)^2+\sum _{j=2}^{t_{i_0,1}^0}(m_{i_0,j}^0)^2+(p_0)^2+\sum _{i \ne i_{\infty}}\sum _{j=1}^{t_{i,1}^{\infty}}(m_{i,j}^{\infty})^2+\sum _{j=2}^{t_{i_{\infty},1}^{\infty}}(m_{i_{\infty},j}^{\infty})^2+(p_{\infty})^2\\
&\ \ +\sum _{i=1}^r\sum _{j=2}^{k_i}(n_j^i)^2+\sum _{i=r+1}^l\sum _{j=1}^{k_i}(n_j^i)^2+\sum _{k=1}^N(p_k)^2-N\{ \dim (mc_{\lambda}(\mathcal{V}))\} ^2\\
&=\sum _{i=1}^{l_0}\sum _{j=1}^{t_{i,1}^0}(m_{i,j}^0)^2-(m_{i_0,1}^0)^2+(p_0)^2+\sum _{i=1}^{l_{\infty}}\sum _{j=1}^{t_{i,1}^{\infty}}(m_{i,j}^{\infty})^2-(m_{i_{\infty},1}^{\infty})^2+(p_{\infty})^2\\
&\ \ +\sum _{i=1}^r\sum _{j=1}^{k_i}(n_j^i)^2-\sum _{i=1}^r(n_1^i)^2+\sum _{k=1}^N(p_k)^2-N\{ \dim (mc_{\lambda}(\mathcal{V}))\} ^2\\
&=\mathop{\rm idx} (E_R)-(m_{i_0,1}^0)^2+(Nm-m_{i_{\infty},1}^{\infty}-d)^2-(m_{i_{\infty},1}^{\infty})^2+(Nm-m_{i_0,1}^0-d)^2-\sum _{i=1}^r(n_1^i)^2\\
&\ \ +\sum _{k=1}^N(Nm-m_{i_0,1}^0-m_{i_{\infty},1}^{\infty}-d+d_k)^2-N\{ (N+1)m-m_{i_0,1}^0-m_{i_{\infty},1}^{\infty}-d\} ^2+Nm^2\\
&=\mathop{\rm idx} (E_R).
\end{align*}
The proof of the theorem has been completed.\ $\square$\\
\\
\textbf{{\Large Acknowledgements}}\\
We would like to express my sincere gratitude to T.Oshima, Y.Haraoka, K.Takemura, D.Yamakawa, K.Hiroe, H.Kawakami and S.Ishizaki for their helpful comments and information about the middle convolution. We wish to thank M.Jimbo, M.Noumi, K.Kajiwara, Y.Ohyama, N.Joshi, T.Masuda, T.Takenawa, T.Tsuda, M.Murata, and Y.Katsushima for discussions and interest. This work is partially supported by JSPS KAKENHI no.24540205.
\bibliographystyle{unsrt}

\end{document}